\documentclass[a4paper,12pt,reqno]{amsart}

\usepackage{amsfonts}
\usepackage{amsmath}
\usepackage{amssymb}
\usepackage{cite}
\usepackage{mathrsfs}
\usepackage{enumitem}
\usepackage[colorlinks]{hyperref}

\numberwithin{equation}{section}

\usepackage{graphicx}

\newtheorem{theorem}{Theorem}[section]
\newtheorem{defi}[theorem]{Definition}
\newtheorem{remark}[theorem]{Remark}

\newtheorem{lemma}[theorem]{Lemma}

\begin{document}
	\allowdisplaybreaks
\title[Discrete time-fractional diffusion equation]
{Time-fractional discrete diffusion equation for Schr\"{o}dinger operator}
\author[A. Dasgupta]{Aparajita Dasgupta}
\address{
	Aparajita Dasgupta:
	\endgraf
	Department of Mathematics
	\endgraf
	Indian Institute of Technology, Delhi, Hauz Khas
	\endgraf
	New Delhi-110016 
	\endgraf
	India
	\endgraf
	{\it E-mail address} {\rm adasgupta@maths.iitd.ac.in}
}
\author[S. S. Mondal]{Shyam Swarup Mondal}
\address{
	Shyam Swarup Mondal:
	\endgraf
	Department of Mathematics
	\endgraf
	Indian Institute of Science, Bengaluru
	\endgraf
	Karnataka-560012 
	\endgraf
	India
	\endgraf
	{\it E-mail address} {\rm mondalshyam055@gmail.com}
}
\author[M. Ruzhansky]{Michael Ruzhansky}
\address{
	Michael Ruzhansky:
	\endgraf
	Department of Mathematics: Analysis, Logic and Discrete Mathematics
	\endgraf
	Ghent University, Belgium
	\endgraf
	and
	\endgraf
	School of Mathematical Sciences
	\endgraf
	Queen Mary University of London
	\endgraf
	United Kingdom
	\endgraf
	{\it E-mail address} {\rm ruzhansky@gmail.com}
}
\author[A. Tushir]{Abhilash Tushir}
\address{
	Abhilash Tushir:
	\endgraf
	Department of Mathematics
	\endgraf
	Indian Institute of Technology, Delhi, Hauz Khas
	\endgraf
	New Delhi-110016 
	\endgraf
	India
	\endgraf
	{\it E-mail address} {\rm abhilash2296@gmail.com}
}

\thanks{The first author was supported by Core Research Grant, RP03890G,  Science and Engineering
	Research Board (SERB), DST,  India. Shyam Swarup Mondal is supported by NBHM post-doctoral fellowship from the Department of Atomic Energy (DAE), Government of India at Indian Institute of Science, Bangalore.  SSM also thanks Ghent Analysis \& PDE Center of Ghent University for the financial support of his visit to Ghent University during which this work has been completed.  M. Ruzhansky is supported by the EPSRC Grants  EP/R003025 by the FWO Odysseus 1 grant G.0H94.18N: Analysis and Partial Differential Equations and by the  Methusalem programme of the Ghent University Special Research Fund (BOF) (Grantnumber
01M01021). The last author is supported by the institute assistantship from the Indian Institute of Technology Delhi, India.}
\date{\today}
\begin{abstract}
 This article aims to investigate the semi-classical analog of the general Caputo-type diffusion equation with time-dependent diffusion coefficient associated with the discrete Schr\"{o}dinger operator, $\mathcal{H}_{\hbar,V}:=-\hbar^{-2}\mathcal{L}_{\hbar}+V$ on the lattice $\hbar\mathbb{Z}^{n},$  where  $V$ is a non-negative multiplication operator and $\mathcal{L}_{\hbar}$ is the discrete Laplacian. 
 %  We prove that the Cauchy problem for the general Caputo-type diffusion equation with regular coefficients is well-posedness in the associated Sobolev-type spaces and very weakly well-posed when distributional coefficients are taken into account.
We establish the well-posedness of the Cauchy problem for the general Caputo-type diffusion equation with a regular coefficient in the associated Sobolev-type spaces. However, it is very weakly well-posed when the diffusion coefficient has a distributional singularity. Finally, we recapture the classical solution  (resp. very weak)  for the general Caputo-type diffusion equation in the semi-classical limit $\hbar\to 0$. 
 \end{abstract}
 
\maketitle
\tableofcontents
\section{Introduction}\label{intro}
The aim of this manuscript is to investigate the semi-classical version of the general Caputo-type diffusion  equation  on  the discrete lattice
\begin{equation*}
	\hbar \mathbb{Z}^{n}=\left\{x \in \mathbb{R}^{n}: x=\hbar k, k \in \mathbb{Z}^{n}\right\},
\end{equation*}
with $\hbar>0$ being a  (small) discretization parameter. Furthermore, we investigate the behaviour of its solutions in the semi-classical limit $\hbar \rightarrow 0$.

The fractional differential operators are nonlocal operators that can be thought of as a generalization of classical differential operators of arbitrary integer orders. These generalized differential equations can be used to model many fields of physics and processes, including fractals, optics, finance, traffic flow, fluid mechanics, mathematical biology, viscoelasticity, control, electrochemistry, signal processing, control theory, and material sciences; see, e.g., \cite{fractal,MAINARDI,highf,torvik84} and the references therein.

Significant attention has been devoted by several researchers to studying the behavior of the solutions of time-fractional heat equations with constant, time-dependent, or space-dependent thermal coefficients in various frameworks due to their widespread application in developing laws of physics while studying real-world phenomena. The modeling of various classical and quantum mechanical issues heavily relies on the nonlocal character of fractional derivatives such as Riemann–Liouville, Caputo, Riesz fractional derivative, etc.   In recent times, there has been a significant increase in interest towards initial value problems for fractional differential equations involving Caputo-type differential operators. We refer to recent papers  \cite{2023niyaz2,sin16}  for the theory of diffusion equations associated with the Caputo-type differential operators.

The general Caputo-type fractional derivative denoted by $\mathbb{D}_{(g)}$ is  defined as
\begin{equation}\label{genfracdef}
	\mathbb{D}_{(g)} u(t)=\frac{d}{d t} \int_0^t g(t-\tau) u(\tau) \mathrm{d} \tau-g(t) u(0)
	,\end{equation}  
where the kernel function $g\in L^{1}_{loc}(0,\infty)$ satisfies  the following conditions:
\begin{enumerate}
	\item\label{c1} the Laplace transform $\widetilde{g}(p)=\mathscr{L}(g ; p)$ of $g$ 
	exists for all $p>0$;
	\item\label{c2} $\widetilde{g}(p)$ is a Stieltjes function;
	\item\label{c3} $\widetilde{g}(p) \rightarrow 0$ and $p \widetilde{g}(p) \rightarrow \infty$ as $p \rightarrow \infty$; and
	\item\label{c4} $\widetilde{g}(p) \rightarrow \infty$ and $p \widetilde{g}(p) \rightarrow 0$ as $p \rightarrow 0$.
\end{enumerate}   
We refer to Section \ref{sec:prelim}  for an extensive discussion of the general Caputo-type fractional derivative. Recently, diffusion equations with the general Caputo-type differential operator have been considered by several authors; see  \cite{chung2018,chung2020,chung2019}. 
Specifically, there is an abundance of literature on the Caputo-type time-fractional diffusion equation with constant, time-dependent, or space-dependent coefficient;  see \cite{fer12,man96,2023niyaz2,2023niyaz,sin16}.  For the time-dependent coefficient settings that differ from the lattice $\mathbb{Z}^n$, see \cite{Dong,Dong2,MR4379761,archil2002,MR3610706,MR3166532,Kubica,MR4469006,Berikbol,vergara,Zacher}.
However, a heat equation with a constant coefficient on the lattice framework $\mathbb{Z}^n$ can be seen in the form of a parabolic Anderson model, see  \cite{MR4076766}. 
%The asymptotic analysis of the solution to the heat equation with constant coefficient on $\mathbb{Z}^n$ is a topic of wide interest due to its application in the field of stochastic analysis.

The primary objective of this paper is to investigate the semi-classical analog of the general Caputo-type diffusion  equation in the Euclidean framework,  which is given by the following Cauchy problem:
\begin{equation}\label{mainpde}
	\left\{\begin{array}{l}
		\mathbb{D}_{(g)} u(t, k)+a(t)\mathcal{H}_{\hbar,V} u(t, k)=f(t, k), \quad(t, k) \in(0, T] \times \hbar \mathbb{Z}^n, \\
		u(0, k)=u_0(k), \quad k \in \hbar \mathbb{Z}^n,
	\end{array}\right.
\end{equation}
where   the diffusion coefficient is $a = a(t) > 0, f$ is the source term,  the operator $\mathbb{D}_{(g)}$ is  given by \eqref{genfracdef}, and the discrete Schr\"{o}dinger operator  $\mathcal{H}_{\hbar,V}$  on  $\hbar \mathbb{Z}^{n}$    is defined by
\begin{equation}\label{dhamil}
	\mathcal{H}_{\hbar,V}u(k):=\left(-\hbar^{-2}\mathcal{L}_{\hbar}+V\right)u(k),\quad k\in\hbar\mathbb{Z}^{n},
\end{equation}
with the discrete Laplacian   $\mathcal{L}_{\hbar}$    given by
\begin{equation}
	\mathcal{L}_{\hbar} u(k):=\sum\limits_{j=1}^{n}\left(u\left(k+\hbar v_{j}\right)+u\left(k-\hbar v_{j}\right)\right)-2 n u(k),\quad k\in\hbar\mathbb{Z}^{n},
\end{equation}
where $v_{j}$ is the $j^{th}$ basis vector in $\mathbb{Z}^{n}$, having all the zeros except for $1$ at the $j^{th}$ component.
Here,  $V$ is a positive multiplication operator by $V(k)$ satisfying the conditions
\begin{equation*}
	V(k)\geq V_{0}>0,\quad \forall~k\in\hbar\mathbb{Z}^{n} \quad \text{and}\quad |V(k)|\to \infty\text{  as  }k\to\infty. 
\end{equation*}
We refer to \cite[Section 3]{Sch:arxiv} for a detailed study on the spectral properties of the discrete Schr\"{o}dinger operator $\mathcal{H}_{\hbar, V}$ on $\hbar \mathbb{Z}^{n}$.    

To the best of our knowledge, the general Caputo-type diffusion equation in the lattice $\hbar\mathbb{Z}^n$ framework has never been considered in the literature so far,  even for the Caputo-type diffusion equation. Our primary objective of the present work is to establish the well-posedness results of the diffusion equation  \eqref{mainpde} in the semi-classical framework of $\hbar \mathbb{Z}^n$.
More precisely, we discuss the well-posedness of the Cauchy problem for the general Caputo-type diffusion equation \eqref{mainpde}  when time-dependent coefficients and source terms are regular (irregular).
Here, we want to emphasize that the discrete wave equation with regular propagation speed and discrete Klein-Gordon with regular potential have already been considered by the authors in \cite{Wave-JDE} and \cite{Klein-Arxiv}, respectively, and prove that both are well-posed in $\ell^2\left(\hbar \mathbb{Z}^n\right)$.     On the other hand,  recently, the authors in \cite{cha23} investigated the well-posedness of the discrete heat equation with regular (irregular) thermal conductivity and discrete tempered distributional data.

In our approach, a special feature is that we also allow the coefficient  $a$ to have distributional irregularities such as $\delta,~ \delta^{2},~ \delta\delta^{\prime}$ or Heaviside function, etc.  Due to Schwartz's impossibility result (see Schwartz \cite{MR64324}),  the singularities of the coefficients would give rise to fundamental mathematical challenges, and in this scenario, the formulation of the Cauchy problem \eqref{mainpde}  may not be feasible from a distributional perspective. Garetto and the third author of this paper proposed a notion of a very weak solution to get around this problem in     \cite{GRweak} and subsequently implemented in    \cite{GR1,GR3,GR2} for various physical models.  

The equation (\ref{mainpde})  is the semi-classical discretization of  the following diffusion equation in the Euclidean framework $\mathbb{R}^n$:
\begin{equation}\label{Eucledian}
	\left\{\begin{array}{l}
		\mathbb{D}_{(g)} u(t,x)+a(t)\mathcal{H}_{V} u(t,x)=f(t, x), \quad (t,x) \in(0, T]\times\mathbb{R}^n, \\
		u(0, x)=u_0(x), \quad x \in \mathbb{R}^{n},
	\end{array}\right.
\end{equation}
where $\mathcal{H}_{V}$ is the usual Schr\"{o}dinger operator on $\mathbb{R}^{n}$ with potential $V$ defined as
\begin{equation}\label{eucd}
	\mathcal{H}_{V}u(x):=\left(-\mathcal{L}+V\right)u(x),\quad x\in\mathbb{R}^{n},
\end{equation}
where $\mathcal{L}$ is the usual Laplacian on $\mathbb{R}^{n}$. 
In the case of a regular coefficient,  using the technique  of    \cite{2023niyaz2} with appropriate modifications (see Section \ref{remark}),  one can   prove that the Cauchy problem  (\ref{Eucledian})    is well-posed in the Sobolev space  $\mathrm{H}^{s}_{\mathcal{H}_{V}}$,
%the space of distributions  $u\in  \mathrm{H}_{\mathcal{H}_{ V}}^{-\infty}$  such that $\left(I+\mathcal{H}_{ V}\right)^{s / 2} u \in L^2 (\mathbb{R}^n)$, and can be  state as: 
%  
defined as 
\begin{equation*}
	\mathrm{H}_{\mathcal{H}_{ V}}^s:=\left\{u \in \mathrm{H}_{\mathcal{H}_{ V}}^{-\infty}:\left(I+\mathcal{H}_{ V}\right)^{s / 2} u \in L^2\left(  \mathbb{R}^n\right)\right\},    \quad s\in \mathbb{R},
\end{equation*}
with the norm $\|f\|_{\mathrm{H}_{\mathcal{H}_{ V}}^s}:=\left\|\left(I+\mathcal{H}_{ V}\right)^{s / 2} f\right\|_{L^2\left(  \mathbb{R}^n\right)}$,  where $\mathrm{H}_{\mathcal{H}_{ V}}^{-\infty}$ is the space of $\mathcal{H}_{ V}$-distributions. The relation between the Sobolev space  associated with the Schr\"{o}dinger operator $\mathrm{H}_{\mathcal{H}_{ V}}^{s}$ and usual Sobolev spaces $\mathrm{H}^{s}(\mathbb{R}^{n})$, particularly for non-negative polynomial potentials, can be understood by using the inequalities obtained by Dziubanski and Glowacki \cite{Dziub}  in the following theorem:
\begin{theorem}\cite{Dziub} Let  $P(x)$ be a nonnegative homogeneous elliptic polynomial on $\mathbb{R}^{n}$ and $V$ is a nonnegative
	polynomial potential.
	For  $1<p<\infty$ and $\alpha>0$ there exist constants $C_1, C_2>0$ such that
	\begin{equation}\label{imp}
		\left\|P(D)^\alpha f\right\|_{L^p}+\left\|V^\alpha f\right\|_{L^p} \leq C_1\left\|(P(D)+V)^{\alpha} f\right\|_{L^p},
	\end{equation}
	and
	\begin{equation*}
		\left\|(P(D)+V)^{\alpha} f\right\|_{L^p} \leq C_2\left\|\left(P(D)^\alpha+V^\alpha\right) f\right\|_{L^p},
	\end{equation*}
	for $f$ in the Schwartz class $\mathcal{S}(\mathbb{R}^{n})$.
\end{theorem}
Now applying the inequality \eqref{imp}  for Schr\"{o}dinger operator $\mathcal{H}_{V}$ with non-negative polynomial potential and using the fact that the Schwartz space $\mathcal{S}(\mathbb{R}^{n})$ is dense in $L^{2}(\mathbb{R}^{n})$, we deduce the following:
\begin{equation}\label{semb}
	\mathrm{H}_{\mathcal{H}_{V}}^{s}(\mathbb{R}^{n})\subseteq \mathrm{H}^{s}(\mathbb{R}^{n}),\quad s>0.
\end{equation}
We refer to \cite{Sch:arxiv}  for a detailed study on the Sobolev spaces for   $\mathcal{H}_{V},\mathrm{H}_{\mathcal{H}_{ V}}^{s}$ and $\mathrm{H}^{-\infty}_{\mathcal{H}_{V}}$. 
We refer to  \cite{nonhm16,nonhm17,ldn1} for Fourier analysis associated with  positive self-adjoint operators with discrete spectrum.

To simplify the notation, we shall write $A\lesssim B$ throughout the article if there
is a constant $C$ that is independent of the  parameters that occur such that $A \leq CB$.
The well-posedness result for the Cauchy problem  \eqref{Eucledian}  with   regular coefficients is given by the following theorem:
\begin{theorem}\label{eucclass}			
	For $s\in\mathbb{R}$, let  $f\in C([0,T];\mathrm{H}^{2+s}_{\mathcal{H}_{V}})$ and $u_{0}\in \mathrm{H}^{2+s}_{\mathcal{H}_{V}}$.  If we suppose that $a\in C([0,T])$  is positive, then the Cauchy problem \eqref{Eucledian} has a unique   solution  $u\in  C([0, T] ; \mathrm{H}_{\mathcal{H}_{V}}^{2+s})$  satisfying $\mathbb{D}_{(g)}u\in  C([0, T] ; \mathrm{H}_{\mathcal{H}_{V}}^{s})$. Furthermore, the solution $u$ also satisfies the estimate
	\begin{equation}
		\|u(t,\cdot)\|_{\mathrm{H}^{2+s}_{\mathcal{H}_{V}}}+	\|\mathbb{D}_{(g)}u(t,\cdot)\|_{\mathrm{H}^{s}_{\mathcal{H}_{V}}}\lesssim\left(1+\|a\|_{C([0,T])}\right)\left[\|u_{0}\|_{\mathrm{H}^{2+s}_{\mathcal{H}_{V}}}+\|f\|_{C([0,T];\mathrm{H}^{2+s}_{\mathcal{H}_{V}})}\right],
	\end{equation}
	for all $t\in[0,T]$.
\end{theorem}
On the other hand, the Cauchy problem \eqref{Eucledian} is very weakly well-posed when the diffusion coefficient has  a distributional singularity.
\begin{theorem}\label{Euwk}
	For $s \in \mathbb{R}$, let $u_{0}  \in\mathrm{H}_{\mathcal{H}_{V}}^{2+s} $. Let  $a$   be  the distribution with support included in $[0, T]$ such that $a \geq a_0>0$ for some positive constant $a_{0}$, and also let the source term
	$f(\cdot, x)$ be a  distribution with support included in $[0, T]$, for all $x\in\mathbb{R}^{n}$.  Then the Cauchy problem \eqref{Eucledian} has a very weak solution  $\left(u_{\varepsilon}\right)_{\varepsilon} \in C([0, T] ; \mathrm{H}_{\mathcal{H}_{V}}^{s})^{(0,1]}$ of order $s$.
\end{theorem}

Concerning the  classical  as well as  very weak solutions to the Cauchy   problem (\ref{mainpde})  on $\hbar \mathbb{Z}^n$, one may ask the following  natural question:
\begin{itemize}
	%\item Can we recapture the classical solution to the heat equation (1.3) when the space variable lies in the Euclidean space R  by allowing  $\hbar\to 0$?
	\item[]	Is it possible to recover the classical solution and very weak solution to the Cauchy problem (\ref{mainpde})  when the spatial variable is taken from the   Euclidean space $\mathbb{R}^n$  by taking   $\hbar\to 0$?
\end{itemize}
In this paper,  we provide affirmative answers to the abovementioned question with some additional Sobolev regularity. 

We conclude the introduction with a brief outline of the organization of the paper. 	
\begin{itemize}
	\item   In Section \ref{sec:results}, we first recall the  basic  Fourier analysis associated with the Schr\"{o}dinger operator $\mathcal{H}_{V}$,  developed by the authors in \cite{Sch:arxiv}. We then present our main results  and their explanations.  The concept of very weak solutions is introduced to address equations that might not have a meaningful solution in the conventional distributional sense.
	
	\item  In Section \ref{sec:prelim}, we recall some essential tools and important results related to the general Caputo-type fractional derivative. 
	
	\item In Section \ref{sec:main results}, we provide proofs of our main results.   In particular,  we establish the wellposedness result to the Cauchy problem \eqref{mainpde} in the case of regular coefficients and also provide the existence, uniqueness, and consistency of the very weak solution to   \eqref{mainpde} in the case of irregular coefficients.

	\item In Section \ref{sec:semiclassic},  we discuss the approximation of the classical solution in the Euclidean framework  $\mathbb{R}^n$  by the solutions in the discrete framework $\hbar \mathbb{Z}^n$. More precisely, 
	we recapture the classical solution and the very weak solution in certain Sobolev-type spaces in the semi-classical limit  $\hbar\to 0$. 
	\item Finally, in Section \ref{remark}, we make a few comments about the classical solution and the very weak solution of the Cauchy problem \eqref{Eucledian}.
\end{itemize}

\section{Main results}\label{sec:results}
In this section, we will summarize our key results for the Cauchy problem \eqref{mainpde} with a regular (irregular) coefficient and a source term. In order to state our main results, some notations and basic preliminaries must be recalled. We refer to Section \ref{sec:prelim} for a thorough exposition of preliminaries. 

The Fourier analysis associated with the Schr\"{o}dinger operator $\mathcal{H}_{\hbar,V}$ was developed by the authors in \cite{Sch:arxiv}. We are simply mentioning a handful of the most essential properties; for detailed  information, see \cite{Sch:arxiv}.  
\begin{lemma}
	Let $V\geq V_{0}>0$ be a multiplication operator satisfying $|V(k)|\to \infty$ as $|k|\to \infty$. Then the discrete Schr\"{o}dinger operator $\mathcal{H}_{\hbar,V}$ has a purely discrete spectrum.
\end{lemma}

% We now proceed to study the Fourier analysis of discrete Hamiltonian operator.
%  In \cite{RS2009,RS2013}, Rabinovich and Roch studied the Fredholm lemmaerty and discrete spectrum of discrete Schr\"{o}dinger opereator (Hamiltonian operator).
Since $-\mathcal{L}_{\hbar}$ is non-negative and $V$ is a positive operator, it follows that $\mathcal{H}_{\hbar,V}$ is a positive operator. Let  $u_{\xi}$ be the eigenfunction of $\mathcal{H}_{\hbar,V}$  associated with the eigenvalue $\lambda_{\xi}$ for each $\xi\in\mathcal{I}_{\hbar}$,  i.e., 
\begin{equation}\label{eigen}
	\mathcal{H}_{\hbar,V}u_{\xi}=\lambda_{\xi}u_{\xi},\quad \text{for all } \xi\in \mathcal{I}_{\hbar},
\end{equation} 
where $\mathcal{I}_{\hbar}$ is a countable set. The collection of eigenvalues, i.e., $\{\lambda_{\xi}> 0:\xi\in\mathcal{I}_{\hbar}\}$  will be the discrete spectrum of the Schr\"{o}dinger operator $\mathcal{H}_{\hbar,V}$. Moreover, we have  
\begin{equation}\label{lambda}
	\Lambda:=\inf\limits_{\hbar>0}\{\lambda_{\xi}:\xi\in\mathcal{I}_{\hbar}\}>0.   
\end{equation}
Consequently,   the  set of eigenfunctions $\{u_{\xi}\}_{\xi\in\mathcal{I}_{\hbar}}$ forms an orthonormal basis for $\ell^{2}(\hbar\mathbb{Z}^{n})$, i.e.,
\begin{equation*}
	\left(u_{\xi},u_{\eta}\right):=\left\{
	\begin{array}{cc}
		1, & \text{if }\xi=\eta, \\
		0, & \text{if }\xi\neq\eta,
	\end{array}\right.
\end{equation*}
where 
\begin{equation*}
	\left(f,g\right):=\sum_{k \in \hbar\mathbb{Z}^{n}}f(k)\overline{g(k)},
\end{equation*}
is the usual inner product of the Hilbert space $\ell^{2}(\hbar\mathbb{Z}^{n})$.

We will now review the spaces of distribution generated by $\mathcal{H}_{\hbar,V}$, its adjoint, and the associated global Fourier analysis. 

The space $\mathrm{H}^{\infty}_{\mathcal{H}_{\hbar,V}}:=\operatorname{Dom}\left(\mathcal{H}_{\hbar,V}^{\infty}\right)$ is called the space of test functions for $\mathcal{H}_{\hbar,V}$,  defined by
\begin{equation*}
	\operatorname{Dom}\left(\mathcal{H}_{\hbar,V}^{\infty}\right):=\bigcap_{k=1}^{\infty} \operatorname{Dom}\left(\mathcal{H}_{\hbar,V}^{k}\right),
\end{equation*}
where $\operatorname{Dom}\left(\mathcal{H}_{\hbar,V}^{k}\right)$ is the domain of the operator $\mathcal{H}_{\hbar,V}^{k}$ defined as
\begin{equation*}
	\operatorname{Dom}\left(\mathcal{H}_{\hbar,V}^{k}\right):=\left\{u \in \ell^{2}(\hbar\mathbb{Z}^{n}): (I+\mathcal{H}_{\hbar,V})^{j} u \in \operatorname{Dom}(\mathcal{H}_{\hbar,V}), j=0,1, \ldots, k-1\right\}.
\end{equation*}
The Fréchet topology of $\mathrm{H}^{\infty}_{\mathcal{H}_{\hbar,V}}$ is given by the family of norms
\begin{equation*}
	\|u\|_{\mathrm{H}^{k}_{\mathcal{H}_{\hbar,V}}}:=\max _{j \leq k}\left\|(I+\mathcal{H}_{\hbar,V})^{j} u\right\|_{\ell^{2}(\hbar\mathbb{Z}^{n})}, \quad k \in \mathbb{N}_{0}, u \in \mathrm{H}^{\infty}_{\mathcal{H}_{\hbar,V}}.
\end{equation*}
The space 
\begin{equation}\label{hvdistr}
	\mathrm{H}^{-\infty}_{\mathcal{H}_{\hbar,V}}:=\mathcal{L}\left(\mathrm{H}^{\infty}_{\mathcal{H}_{\hbar,V}}, \mathbb{C}\right),	
\end{equation}
of all continuous linear functionals on $\mathrm{H}^{\infty}_{\mathcal{H}_{\hbar,V}}$ is called the space of $\mathcal{H}_{\hbar,V}$-distributions. For $w \in \mathrm{H}^{-\infty}_{\mathcal{H}_{\hbar,V}}$ and $u \in \mathrm{H}^{\infty}_{\mathcal{H}_{\hbar,V}}$, we shall write
\begin{equation*}
	w(u)=( w, u) =\sum\limits_{k\in\hbar\mathbb{Z}^{n}}w(k)\overline{u(k)}.
\end{equation*}
For any $u \in \mathrm{H}^{\infty}_{\mathcal{H}_{\hbar,V}}$, the functional
\begin{equation*}
	\mathrm{H}^{\infty}_{\mathcal{H}_{\hbar,V}} \ni v \mapsto (v,u),
\end{equation*}
is a $\mathcal{H}_{\hbar,V}$-distribution, which gives an embedding $u \in \mathrm{H}^{\infty}_{\mathcal{H}_{\hbar,V}} \hookrightarrow \mathrm{H}^{-\infty}_{\mathcal{H}_{\hbar,V}}.$ 
\begin{defi}[\textbf{Schwartz Space} $\mathcal{S}(\mathcal{I}_{\hbar})$]
	Let $\mathcal{S}(\mathcal{I}_{\hbar})$ denote the space of rapidly decaying functions $\varphi:\mathcal{I}_{\hbar}\to \mathbb{C}$, i.e.,  $\varphi \in \mathcal{S}(\mathcal{I}_{\hbar})$ if for any $M<\infty$, there exists a constant $C_{\varphi, M}$ such that
	\begin{equation*}
		|\varphi(\xi)| \leq C_{\varphi, M}\langle\xi\rangle^{-M}, \quad\text{for all } \xi\in\mathcal{I}_{\hbar},
	\end{equation*}
	where we denote
	\begin{equation*}
		\langle\xi\rangle:=\left(1+\lambda_{\xi}\right)^{\frac{1}{2 }}.
	\end{equation*}
\end{defi}
The topology on $\mathcal{S}(\mathcal{I}_{\hbar})$ is given by the family of seminorms $p_{k}$, where $k \in \mathbb{N}_{0}$ and
\begin{equation*}
	p_{k}(\varphi):=\sup _{\xi \in \mathcal{I}_{\hbar}}\langle\xi\rangle^{k}|\varphi(\xi)| .
\end{equation*}
We now define the $\mathcal{H}_{\hbar,V}$-Fourier transform on $\mathrm{H}^{\infty}_{\mathcal{H}_{\hbar,V}}$.
\begin{defi}
	The $\mathcal{H}_{\hbar,V}$-Fourier transform  $	\mathcal{F}_{\mathcal{H}_{\hbar,V}}: \mathrm{H}^{\infty}_{\mathcal{H}_{\hbar,V}} \rightarrow \mathcal{S}(\mathcal{I}_{\hbar})$ is given by the formula
	\begin{equation*}
		\left(\mathcal{F}_{\mathcal{H}_{\hbar,V}} f\right)(\xi)=	\widehat{f}(\xi):=\sum\limits_{k\in\hbar\mathbb{Z}^{n}} f(k) \overline{u_{\xi}(k)}, \quad \xi\in\mathcal{I}_{\hbar},
	\end{equation*}
	where $u_{\xi}$ satisfies \eqref{eigen}. It can also be extended to $\mathrm{H}^{-\infty}_{\mathcal{H}_{\hbar,V}}$ in the usual way.
\end{defi}
Moreover, the $\mathcal{H}_{\hbar,V}$-Fourier transform $\mathcal{F}_{\mathcal{H}_{\hbar,V}}:\mathrm{H}^{\infty}_{\mathcal{H}_{\hbar,V}} \rightarrow \mathcal{S}(\mathcal{I}_{\hbar})$ is a homeomorphism and its  inverse $
\mathcal{F}_{\mathcal{H}_{\hbar,V}}^{-1}:\mathcal{S}\left(\mathcal{I}_{\hbar}\right)\rightarrow	\mathrm{H}^{\infty}_{\mathcal{H}_{\hbar,V}},$
is given by 
\begin{equation*}
	\left(\mathcal{F}_{\mathcal{H}_{\hbar,V}}^{-1}g\right)(k):=\sum\limits_{\xi\in\mathcal{I}_{\hbar}}g(\xi)u_{\xi}(k),\quad g\in \mathcal{S}(\mathcal{I}_{\hbar}),
\end{equation*}
so that the $\mathcal{H}_{\hbar,V}$-Fourier inversion formula becomes
\begin{equation*}
	f(k)=\sum_{\xi \in \mathcal{I}_{\hbar}} \widehat{f}(\xi) u_{\xi}(k), \quad \text { for all } f \in \mathrm{H}^{\infty}_{\mathcal{H}_{\hbar,V}}.
\end{equation*}
The $\mathcal{H}_{\hbar,V}$-Plancherel formula takes the form
\begin{equation}\label{planch}
	\sum_{k \in \hbar\mathbb{Z}^{n}}|f(k)|^{2}=\sum_{\xi \in \mathcal{I}_{\hbar}}|\widehat{f}(\xi)|^{2}.
\end{equation}
The  $\mathcal{H}_{\hbar,V}$-Fourier transform of the discrete Schr\"{o}dinger operator $\mathcal{H}_{\hbar,V}$ is
\begin{equation}\label{symbol}
	\widehat{\mathcal{H}_{\hbar,V}f}(\xi)=(\mathcal{H}_{\hbar,V}f,u_{\xi})=(f,\mathcal{H}_{\hbar,V}u_{\xi})=(f,\lambda_{\xi}u_{\xi})=\lambda_{\xi}\widehat{f}(\xi),\quad \xi\in\mathcal{I}_{\hbar},
\end{equation} 
since $\mathcal{H}_{\hbar,V}$ is a self-adjoint operator.  Furthermore,
for $s \in \mathbb{R}$, we introduce the Sobolev space $\mathrm{H}_{\mathcal{H}_{\hbar, V}}^s$ defined by
\begin{equation}\label{sobdef}
	\mathrm{H}_{\mathcal{H}_{\hbar, V}}^s:=\left\{u \in \mathrm{H}_{\mathcal{H}_{\hbar, V}}^{-\infty}:\left(I+\mathcal{H}_{\hbar, V}\right)^{s / 2} u \in \ell^2\left(\hbar \mathbb{Z}^n\right)\right\},    
\end{equation}
with the norm $\|f\|_{\mathrm{H}_{\mathcal{H}_{\hbar, V}}^s}:=\left\|\left(I+\mathcal{H}_{\hbar, V}\right)^{s / 2} f\right\|_{\ell^2\left(\hbar \mathbb{Z}^n\right)}$,  where $\mathrm{H}_{\mathcal{H}_{\hbar, V}}^{-\infty}$ is the space of $\mathcal{H}_{\hbar,V}$-distributions.
Combining the Plancherel identity \eqref{planch} with \eqref{symbol}, we obtain
\begin{equation}\label{fnorm}
	\|f\|_{\mathrm{H}_{\mathcal{H}_{\hbar,V}}^{s}}:=\left\|(I+\mathcal{H}_{\hbar,V})^{s / 2} f\right\|_{\ell^{2}(\hbar\mathbb{Z}^{n})}=\left(\sum\limits_{\xi\in\mathcal{I}_{\hbar}}\left(1+\lambda_{\xi}\right)^{s}|\widehat{f}(\xi)|^{2}\right)^{\frac{1}{2}}.
\end{equation}

For more details about the discrete Schr\"{o}dinger operator and the associated Sobolev spaces, one can refer to \cite[Section 3]{Sch:arxiv}. 

Let us begin by outlining the  well-posedness result for the Cauchy problem \eqref{mainpde} with
regular coefficient $a$ and the source term $f$:
\begin{theorem}\label{mthm1}
	For $s\in\mathbb{R}$, let  $f\in C([0,T];\mathrm{H}^{2+s}_{\mathcal{H}_{\hbar,V}})$ and $u_{0}\in \mathrm{H}^{2+s}_{\mathcal{H}_{\hbar,V}}$.  If we suppose that $a\in C([0,T])$  is positive, then the Cauchy problem \eqref{mainpde} has a unique  solution  $u\in C([0,T];\mathrm{H}^{2+s}_{\mathcal{H}_{\hbar,V}})$ satisfying $\mathbb{D}_{(g)}u\in C([0,T];\mathrm{H}^{s}_{\mathcal{H}_{\hbar,V}})$. Furthermore, the solution $u$  can be represented by 
	\begin{equation}\label{specrep}
		u(t,k)=\sum\limits_{\xi\in\mathcal{I}_{\hbar}}\widehat{u}(t,\xi)u_{\xi}(k), \quad (t,k)\in [0,T]\times \hbar\mathbb{Z}^{n},
	\end{equation}
	where $u_{\xi}$ are the eigenfunctions of $\mathcal{H}_{\hbar,V}$  given in \eqref{eigen}, and also satisfies the estimate
	\begin{equation}\label{reqesttt}
		\|u(t,\cdot)\|_{\mathrm{H}^{2+s}_{\mathcal{H}_{\hbar,V}}}+	\|\mathbb{D}_{(g)}u(t,\cdot)\|_{\mathrm{H}^{s}_{\mathcal{H}_{\hbar,V}}}\leq C_{a}\left[\|u_{0}\|_{\mathrm{H}^{2+s}_{\mathcal{H}_{\hbar,V}}}+\|f\|_{C([0,T];\mathrm{H}^{2+s}_{\mathcal{H}_{\hbar,V}})}\right],
	\end{equation}
	for all $t\in[0,T]$, where the constant $C_{a}$ is given by
	\begin{equation}\label{cnstthm}
		C_{a}=(1+\|a\|_{C([0,T])})\max\left\{1,\frac{1}{\Lambda a_{0}}\right\},   
	\end{equation}
	with $\Lambda:=\inf\limits_{\hbar>0}\{\lambda_{\xi}:\xi\in\mathcal{I}_{\hbar}\}>0$.
\end{theorem}
Next, we will permit singularities for the coefficient $a$ and the source term $f$ in the time variable.
Following the technique from \cite{GRweak}, by regularizing the distributional coefficient $a$ in the manner described below, we can form a net of smooth functions $(a_{\varepsilon})_{\varepsilon}$ given by:
\begin{equation*}
	a_{\varepsilon}(t):=(a*\psi_{\omega(\varepsilon)})(t),\quad  \psi_{\varepsilon}(t)=\frac{1}{\omega(\varepsilon)}\psi\left(\frac{t}{\omega(\varepsilon)}\right), \quad \varepsilon\in(0,1],
\end{equation*}
where $\omega(\varepsilon)\to 0$ as $\omega\to 0$, and $\psi$ is a Friedrichs-mollifier, i.e., $\psi\in C_{0}^{\infty}(\mathbb{R}),\psi\geq 0$, and $\int_{\mathbb{R}}\psi=1$.

For a net of functions/distributions, the notions of moderateness and negligibility are as follows:
\begin{defi}[Moderateness]
	(i)	A net $\left(a_{\varepsilon}\right)_{\varepsilon} \in C(\mathbb{R})^{(0,1]}$  is said to be $C$-moderate if, for all  $K \Subset \mathbb{R}$, there exists $N \in \mathbb{N}_{0}$  such that
	\begin{equation*}	\left\| a_{\varepsilon}\right\|_{C(K)} \lesssim \varepsilon^{-N},
	\end{equation*}
	for all $\varepsilon \in(0,1]$.\\
	(ii)	A net $\left(u_{\varepsilon}\right)_{\varepsilon} \in C([0,T];\mathrm{H}_{\mathcal{H}_{\hbar,V}}^{s})^{(0,1]}$ is said to be $C([0, T] ; \mathrm{H}_{\mathcal{H}_{\hbar,V}}^{s})$-moderate if   there exists $N \in \mathbb{N}_{0}$ such that
	\begin{equation*}
		\|u_{\varepsilon}\|_{C([0,T];\mathrm{H}_{\mathcal{H}_{\hbar,V}}^{s})} \lesssim \varepsilon^{-N},
	\end{equation*}
	for all $\varepsilon \in(0,1]$.
	
\end{defi}
\begin{defi}[Negligibility]
	(i) 	A net $\left(a_{\varepsilon}\right)_{\varepsilon} \in C(\mathbb{R})^{(0,1]}$ is said to be $C$-negligible if, for all $K \Subset \mathbb{R}$ and   $q\in\mathbb{N}_{0}$, the following condition is satisfied
	\begin{equation*}
		\left\| a_{\varepsilon}\right\|_{C(K)} \lesssim \varepsilon^{q},
	\end{equation*}
	for all $\varepsilon \in(0,1]$.\\
	(ii) A net $\left(u_{\varepsilon}\right)_{\varepsilon} \in C([0, T] ; \mathrm{H}_{\mathcal{H}_{\hbar,V}}^{s})^{(0,1]}$ is said to be $C([0, T] ; \mathrm{H}_{\mathcal{H}_{\hbar,V}}^{s})$-negligible if for all $q\in\mathbb{N}_{0}$, the following condition is satisfied 
	\begin{equation*}
		\|u_{\varepsilon}\|_{C([0,T];\mathrm{H}_{\mathcal{H}_{\hbar,V}}^{s})} \lesssim \varepsilon^{q},
	\end{equation*}
	for all  $\varepsilon \in(0,1]$.
\end{defi}
\noindent
The constant appearing for $``\lesssim"$ in the definition of negligibility may depend on $q$ but not on $\varepsilon$. Furthermore, using the structure theorems for the distributions, we have
\begin{equation*}\label{embedding}
	\text{compactly supported distributions } \mathcal{E}^{\prime}(\mathbb{R})\subset \{C\text{-moderate families}\}.   
\end{equation*}
Therefore, even though a solution to a Cauchy problem might not exist in the space on the left side of the above embedding, it might exist in the space on the right side. 

We may now introduce the concept of a very weak solution in our settings in the following manner:
\begin{defi}\label{vwkdef}
	Let $s \in \mathbb{R}$ and $u_{0} \in \mathrm{H}_{\mathcal{H}_{\hbar,V}}^{2+s}.$ The net $\left(u_{\varepsilon}\right)_{\varepsilon} \in$ $C([0, T] ; \mathrm{H}_{\mathcal{H}_{\hbar,V}}^{2+s})^{(0,1]}$ is said to be a very weak solution of order $s$ of the Cauchy problem \eqref{mainpde}, if there exist
	\begin{enumerate}
		\item $C$-moderate regularization $a_{\varepsilon}$    of the coefficient $a$; and
		\item $C([0, T] ; \mathrm{H}_{\mathcal{H}_{\hbar,V}}^{2+s})$-moderate regularization $f_{\varepsilon}$ of the source term $f$,
	\end{enumerate}
	such that $\left(u_{\varepsilon}\right)_{\varepsilon}$ is $C([0, T] ; \mathrm{H}_{\mathcal{H}_{\hbar,V}}^{2+s})$-moderate and solves the regularized problem
	\begin{equation}\label{reg}
		\left\{\begin{array}{l}
			\mathbb{D}_{(g)} u_{\varepsilon}(t, k)+a_{\varepsilon}(t)\mathcal{H}_{\hbar,V}u_{\varepsilon}(t, k)=f_{\varepsilon}(t, k),\quad(t, k) \in(0, T] \times \hbar\mathbb{Z}^{n}, \\
			u_{\varepsilon}(0, k)=u_{0}(k),\quad k \in \hbar\mathbb{Z}^{n},
		\end{array}\right.
	\end{equation}
	for all $\varepsilon \in(0,1]$.
\end{defi}
Noting that the distributions $a$ and $ b$, are distributions with domain $[0, T]$,  it suffices to take into account that $\operatorname{supp}(\psi) \subseteq K,$ where $K = [0, T]$, throughout the article.
A distribution $a$ is said to be  positive if 
\begin{equation}
	\langle a, \psi\rangle \geq 0~ \text{for all}~\psi\in C^{\infty}_{0}(\mathbb{R}), \psi\geq 0,
\end{equation}
where $\langle a, \psi\rangle$
denotes the action of the distribution $a$
on the test function 
$\psi$.  In a similar way, a distribution $a$ is said to be strictly positive if there exists a  constant $\alpha>0$ such that   $a\geq\alpha>0$, where the meaning of the expression $a\geq \alpha$ is that the distribution $a-\alpha$ is positive, i.e., 
\begin{equation*}
	\langle a-\alpha, \psi\rangle \geq 0, \quad \text{for all } \psi \in C_0^{\infty}(\mathbb{R}),\psi\geq 0.
\end{equation*}
Furthermore, the following definition may be employed  to understand the uniqueness of a very weak solution to the Cauchy problem \eqref{mainpde}:
\begin{defi}\label{uniquedef}
	We say that the Cauchy problem \eqref{mainpde} has a $C([0, T] ; \mathrm{H}_{\mathcal{H}_{\hbar,V}}^{2+s})$-unique very weak solution, if
	\begin{enumerate}
		\item for all  $C$-moderate nets $a_{\varepsilon},\tilde{a}_{\varepsilon}$ such that $(a_{\varepsilon}-\tilde{a}_{\varepsilon})_{\varepsilon}$ is $C$-negligible; and
		\item for all $C([0, T] ; \mathrm{H}_{\mathcal{H}_{\hbar,V}}^{2+s})$-moderate nets $f_{\varepsilon},\tilde{f}_{\varepsilon}$ such that $(f_{\varepsilon}-\tilde{f}_{\varepsilon})_{\varepsilon}$ is \\ $C([0, T] ; \mathrm{H}_{\mathcal{H}_{\hbar,V}}^{2+s})$-negligible,
	\end{enumerate}
	the net $(u_{\varepsilon}-\tilde{u}_{\varepsilon})_{\varepsilon}$ is $C([0, T] ; \mathrm{H}_{\mathcal{H}_{\hbar,V}}^{2+s})$-negligible, where
	$(u_{\varepsilon})_{\varepsilon}$ and $(\tilde{u}_{\varepsilon})_{\varepsilon}$ 	 are the families of
	solutions corresponding to the 	$\varepsilon$-parametrised problems
	\begin{equation}\label{reg1}
		\left\{\begin{array}{l}
			\mathbb{D}_{(g)} u_{\varepsilon}(t, k)+a_{\varepsilon}(t)\mathcal{H}_{\hbar,V}u_{\varepsilon}(t, k)=f_{\varepsilon}(t, k),\quad(t, k) \in(0, T] \times \hbar\mathbb{Z}^{n}, \\
			u_{\varepsilon}(0, k)=u_{0}(k),\quad k \in \hbar\mathbb{Z}^{n},
		\end{array}\right.
	\end{equation}
	and
	\begin{equation}\label{reg2}
		\left\{\begin{array}{l}
			\mathbb{D}_{(g)} \tilde{u}_{\varepsilon}(t, k)+\tilde{a}_{\varepsilon}(t)\mathcal{H}_{\hbar,V}\tilde{u}_{\varepsilon}(t, k)=\tilde{f}_{\varepsilon}(t, k),\quad(t, k) \in(0, T] \times \hbar\mathbb{Z}^{n}, \\
			\tilde{u}_{\varepsilon}(0, k)=u_{0}(k),\quad k \in \hbar\mathbb{Z}^{n},
		\end{array}\right.
	\end{equation}
	respectively.
\end{defi}
The existence and uniqueness of a very weak solution for the
Cauchy problem \eqref{mainpde} with the distributional coefficient as well as the source term can be stated as follows:
\begin{theorem}[Existence and uniqueness]\label{ext}
	For $s\in\mathbb{R}$, let  $u_{0} \in \mathrm{H}_{\mathcal{H}_{\hbar,V}}^{2+s}$.	Let   $a$  be a   distribution with support contained in $[0, T]$ such that $a\geq a_{0}>0$ for some positive constant $a_{0}$, and also let  the source term $f(\cdot,k)$ be a distribution with  support contained in $[0,T]$, for all $k\in\hbar\mathbb{Z}^{n}$. Then the Cauchy problem \eqref{mainpde} has a $C([0,T];\mathrm{H}_{\mathcal{H}_{\hbar,V}}^{2+s})$-unique very weak solution of order $s$.
\end{theorem}
Our next theorem is devoted to the consistency of the very weak solution with the classical solution, i.e., whenever the classical solution $u$ of the Cauchy problem \eqref{mainpde} exists, the very weak solution $(u_{\varepsilon})_{\varepsilon}\to u$ in $C([0,T];\mathrm{H}_{\mathcal{H}_{\hbar,V}}^{2+s})$ as $\varepsilon\to 0$.
\begin{theorem}[Consistency]\label{cnst}
	For $s\in \mathbb{R}$, let  $f \in C([0, T], \mathrm{H}_{\mathcal{H}_{\hbar,V}}^{2+s})$ and $u_{0} \in \mathrm{H}_{\mathcal{H}_{\hbar,V}}^{2+s}$. Let $(u_{\varepsilon})_{\varepsilon}$ be the very weak solution and $u$ be the classical solution of the Cauchy problem \eqref{mainpde} given by Theorem \ref{mthm1}. If we suppose that $a \in C([0, T])$ satisfies $\inf\limits_{t \in[0, T]}a(t)= a_{0}>0$, then for any regularization $(a_{\varepsilon})_{\varepsilon}$ and $ (f_{\varepsilon})_{\varepsilon}$, we have
	\begin{equation*}
		u_{\varepsilon} \to u \text{ in }C([0, T]; \mathrm{H}_{\mathcal{H}_{\hbar,V}}^{2+s}) ,\quad \varepsilon\to 0.
	\end{equation*}
\end{theorem}
In the following theorem,   we approximate the classical solution in Euclidean settings by the solutions in discrete settings, namely,  we  prove that under the assumption that the solutions in $\mathbb{R}^{n}$ exist, they can be recovered in the limit as $\hbar\to 0$. 
\begin{theorem}\label{semlimit}
	Let $V\geq V_{0}>0$ be a positive polynomial potential in \eqref{dhamil}.	Let $u$ and $v$ be the solutions for  the  Cauchy problems \eqref{mainpde} on $\hbar \mathbb{Z}^{n}$ and \eqref{Eucledian} on $\mathbb{R}^{n}$, respectively, with the same Cauchy data  $u_{0}$ and the source term $f$. Assume  that the initial Cauchy data $u_{0} \in \mathrm{H}_{\mathcal{H}_{V}}^{2+s} $  with $s>2+\frac{n}{2}$  satisfies $u_{0}^{(4v_{j})}\in \mathrm{H}_{\mathcal{H}_{V}}^{2+s} $ for all $j=1,\dots,n$, where $v_{j}$'s are the $j^{th}$ standard basis vectors in $\mathbb{Z}^{n}$.  Then for every $t \in[0, T]$, we have
	\begin{equation*}
		\|v(t)-u(t )\|_{\mathrm{H}_{\mathcal{H}_{\hbar,V}}^{2+s}}+\left\| \mathbb{D}_{(g)}v(t)-\mathbb{D}_{(g)} u(t)\right\|_{\mathrm{H}_{\mathcal{H}_{\hbar,V}}^{s}} \rightarrow 0, \text { as } \hbar \rightarrow 0,
	\end{equation*}
	and the convergence is uniform on $[0,T]$.
\end{theorem}
\begin{remark}
	The evaluation of initial data and source term from \eqref{Eucledian} on the lattice $\hbar\mathbb{Z}^{n}$ is the initial Cauchy data $u_{0}$ and the source term $f$ of the Cauchy problem \eqref{mainpde}.
\end{remark}
The analogous of the above semi-classical theorem for a very weak solution can be stated as follows:
\begin{theorem}\label{vvyksemlimit}
	Let $V\geq V_{0}>0$ be a positive polynomial potential in \eqref{dhamil}. Let $(u_{\varepsilon})_{\varepsilon}$ and $(v_{\varepsilon})_{\varepsilon}$ be the very weak solutions of the Cauchy problems \eqref{mainpde} on $\hbar \mathbb{Z}^{n}$ and \eqref{Eucledian} on $\mathbb{R}^{n}$, respectively, with the same Cauchy data $u_{0}$ and the source term $f$. Assume the initial Cauchy data $u_{0}\in \mathrm{H}_{\mathcal{H}_{V}}^{2+s}$  with $s>2+\frac{n}{2}$ satisfies $u_{0}^{(4v_{j})}\in \mathrm{H}_{\mathcal{H}_{V}}^{2+s}$ for all $j=1,\dots,n$, where $v_{j}$'s are the $j^{th}$ standard basis vectors in $\mathbb{Z}^{n}$.  Then for every $\varepsilon\in(0,1]$ and $t\in[0,T]$, we have
	\begin{equation*}
		\|v_{\varepsilon}(t)-u_{\varepsilon}(t)\|_{\mathrm{H}_{\mathcal{H}_{\hbar,V}}^{2+s}}+\left\|\mathbb{D}_{(g)}v_{\varepsilon}(t)-\mathbb{D}_{(g)} u_{\varepsilon}(t)\right\|_{\mathrm{H}_{\mathcal{H}_{\hbar,V}}^{s}} \rightarrow 0, \text { as } \hbar \rightarrow 0,
	\end{equation*}
	where the convergence is uniform on $[0,T]$ but pointwise for $\varepsilon\in(0,1]$.
\end{theorem}
\section{Preliminaries}\label{sec:prelim}
In this section, we will recall some necessary tools and results related to the general Caputo-type fractional derivative and discrete Schr\"{o}dinger operator that will be necessary for the analysis that we will follow.
\begin{defi}\cite{kilbas}
	The generalized Mittag-Leffler function denoted by $E_{\alpha,\beta}(z)$ is defined by
	\begin{equation*}
		E_{\alpha,\beta}(z)=\sum_{k=0}^{\infty}\frac{z^{k}}{\Gamma(\alpha k+\beta)},\quad (z,\beta\in \mathbb{C}, \mathrm{Re}(\alpha)>0).
	\end{equation*}
\end{defi}
\begin{defi}\cite{Bern}
	A real-valued function $f$ on $(0,\infty)$ is said to be completely monotone if $f$ is infinitely differentiable   and satisfies
	\begin{equation*}
		(-1)^{n}f^{(n)}(t)\geq 0, \quad \text{for all } n\in \mathbb{N}_{0} \text{ and }t>0.
	\end{equation*}
\end{defi}
\noindent Here are some examples of completely monotone functions on $(0,\infty)$:
\begin{equation*}
	c, e^{-ct},\frac{1}{2\sqrt{t}},\frac{1}{t+c^{2}},\quad c>0.
\end{equation*}
\begin{defi}\cite{Bern}
	The Laplace transform of a non-negative function $f$ on $[0,\infty)$ or a measure $\mu$ on the line $[0,\infty)$ is defined as
	\begin{equation*}
		\mathscr{L}(f ; p):=\int_0^{\infty} e^{-p \tau} f(\tau) \mathrm{d} \tau \quad \text { or } \quad \mathscr{L}(\mu ; p):=\int_{0}^{\infty} e^{-p \tau} \mu(\mathrm{d} \tau),\quad p>0,
	\end{equation*}
	respectively, wherever these integrals converge.
\end{defi}
\subsection{Generalized Caputo-type fractional derivative} The general Caputo-type fractional derivative was introduced by Kochubei in \cite{2011Koch} and subsequently considered by many authors for relaxation equations, diffusion equations, evolution equations, etc. (see \cite{chung2018,chung2020,chung2019}). We must recall the Stieltjes function in order to define the general Caputo-type fractional derivative.
\begin{defi}\label{steilt}
	A non-negative function $f$ on $(0,\infty)$ is said to be a Stieltjes function if $f$ admits the integral representation
	\begin{equation*}
		f(t)=\frac{a}{t}+b+\int_{0}^{\infty} \frac{1}{t+\tau} \mu(\mathrm{d} \tau),    
	\end{equation*}
	where $a, b \geq 0$ are non-negative constants and $\mu$ is a positive Borel measure on $(0, \infty)$ such that 
	\begin{equation*}
		\int_{0}^{\infty}\frac{1}{1+\tau} \mu(\mathrm{d} \tau)<\infty.   
	\end{equation*}
\end{defi}
\noindent Here are  some examples of Stieltjes functions on $(0,\infty)$: 
\begin{equation*}
	c,~\frac{1}{t+c},~\frac{1+c}{t+c},~\frac{1}{t}\log(1+t),~\frac{1}{\sqrt{t}}\arctan(t),\quad c>0. 
\end{equation*}
The class of Stieltjes functions forms a subclass of completely monotone functions. Furthermore, in view of \cite[Proposition 1.2]{Bern}, the constants $a,b$ and the measure $\mu$ in the Definition \ref{steilt} are uniquely determined by $f$. For more information about the Stieltjes function, one can refer to \cite{Bern}. Now we are in a position to define the general Caputo-type fractional derivative. 
\begin{defi}\cite{2011Koch}\label{Caputo-type}
	The general Caputo-type fractional derivative denoted by $\mathbb{D}_{(g)}$ is defined as
	\begin{equation*}
		\mathbb{D}_{(g)} u(t)=\frac{d}{d t} \int_0^t g(t-\tau) u(\tau) \mathrm{d} \tau-g(t) u(0)=\int_{0}^{t}g(t-\tau)u^{\prime}(\tau)\mathrm{d}\tau
		,\end{equation*}  
	where the kernel function $g\in L^{1}_{loc}(0,\infty)$ satisfies the conditions \eqref{c1}-\eqref{c4} in the introduction.
\end{defi}
\noindent 
Depending on the choice of a kernel function $g$, there are several variants of the Caputo-type fractional derivative; we list a few of them below.  For $0<\alpha<1$, the general Caputo-type fractional derivative is called:
\begin{enumerate}
	\item\label{CDD} Caputo-Dzhrbashyan derivative when $g(t)=\dfrac{t^{-\alpha}}{\Gamma(1-\alpha)}$;
	\item Caputo-Fabrizio derivative when $g(t)=\dfrac{M(\alpha)}{1-\alpha}\exp\left(-\dfrac{\alpha}{1-\alpha}t\right)$; and
	\item Atangana-Baleanu derivative when $g(t)=\dfrac{B(\alpha)}{1-\alpha}E_{\alpha}\left(-\dfrac{\alpha}{1-\alpha}t\right)$,
\end{enumerate}
where $M(\alpha)$ and $B(\alpha)$ are normalization factors satisfying $M(0)=M(1)=B(0)=B(1)=1$, see \cite{diet20}. For more information about the general Caputo-type fractional derivative, one can refer to \cite{2011Koch}.

Next, we will review some results for the Cauchy problem associated with the general Caputo-type fractional derivative $\mathbb{D}_{(g)}$:
\begin{equation}\label{homchy}
	\left\{\begin{array}{l}
		\mathbb{D}_{(g)} w(t)+\lambda w(t)=f(t), \quad t>0,\lambda>0, \\
		w(0)=w_{0},
	\end{array}\right.
\end{equation}
and make some appropriate remarks. The following lemma is dedicated to the corresponding homogeneous fractional differential equation.
\begin{lemma}\label{homthm}\cite[Theorem 2]{2011Koch}  The fractional differential equation \eqref{homchy} with initial Cauchy data $w_{0}=1$ and the source term $f\equiv 0$, 
	has a  unique solution $w_{\lambda}(t)$ satisfying the following:
	\begin{enumerate}
		\item\label{P1} $w_{\lambda}$ is continuous on $[0,\infty)$;
		\item\label{P2} $w_{\lambda}$ is infinitely differentiable on $(0,\infty)$;  and
		\item\label{P3} $w_{\lambda}$ is completely monotone on $(0,\infty)$.
	\end{enumerate}
\end{lemma}
In the case of the Caputo-Dzhrbashyan derivative, the generalized Mittag-Leffler function can be utilized to express the solution $w_{\lambda}$ as follows: 
\begin{equation*}
	w_{\lambda}(t)=E_{\alpha,1}(-\lambda t^{\alpha}),\quad t>0.  
\end{equation*}
\begin{remark}\label{firstremark}
	In continuation of the above theorem, using the monotonicity of $w_{\lambda}$ and the Karamata-Feller Tauberian theorem in \cite{fel:book}, the author also proved that 
	\begin{equation*}
		w_{\lambda}(t)\to 1,\quad \text{as }t\to 0^{+}.
	\end{equation*}
\end{remark}
\begin{remark}\label{2ndremark}
	In light of Lemma \ref{homthm} and the above remark, it is easy to note that the solution $w_{\lambda}(t)$ satisfies
	\begin{equation*}
		0\leq w_{\lambda}(t)\leq 1, \quad \text{ for all }t>0.
	\end{equation*}
\end{remark}
\begin{lemma}\cite[Lemma 3.1]{chung2020}\label{nonthm}
	Let $\lambda>0, T>0$, and $f \in C([0, T])$. Then the fractional differential equation \eqref{homchy} with the initial Cauchy data $w_{0}=0$,
	has a unique solution in $C([0, T])$. In particular, the solution has the form
	\begin{equation*}
		w(t)=-\frac{1}{\lambda} \int_0^t  f(s)w_{\lambda}^{\prime}(t-s) \mathrm{d} s,
	\end{equation*}
	where $w_{\lambda}(t)$ is the solution of the homogeneous Cauchy problem obtained in Lemma \ref{homthm}.
\end{lemma}
Similar to Lemma \ref{homthm}, the solution $w(t)$ obtained in the above lemma in the case of the Caputo-Dzhrbashyan derivative can be expressed as
\begin{equation*}
	w(t)=-\frac{1}{\lambda}\int_{0}^{t}f(s)\left(\frac{d}{ds}E_{\alpha,1}(-\lambda(t-s)^{\alpha})\right)\mathrm{d}s.
\end{equation*}
\begin{remark}\label{3rdremark}
	Combining the Lemmas \ref{homthm} and \ref{nonthm}, the fractional differential equation \eqref{homchy} has a unique solution in $C([0,T])$ expressed as
	\begin{equation*}
		w(t)=w_{0}w_{\lambda}(t)-\frac{1}{\lambda} \int_0^t  f(s)w_{\lambda}^{\prime}(t-s) \mathrm{d} s.
	\end{equation*}
\end{remark}
\begin{remark}\label{4thremark}
	Using the fact that $w_{\lambda}^{\prime}(t)\leq 0$ for all $t\in[0,T]$, we can conclude that the solution $w(t)\geq(\leq)0$ for all $t\in[0,T]$, if the initial Cauchy data $w_{0}\geq(\leq)0$ and the source term $f\geq(\leq)0$ for all $t\in[0,T]$.
\end{remark}
\section{Proof of the main results}\label{sec:main results}
In this section,  we first prove the existence and uniqueness of the classical solution $u$ of the Cauchy problem \eqref{mainpde}, and then we will conclude the proof of Theorem   \ref{mthm1} by establishing the regularity estimates. Throughout this section, we will specify the coefficients as $u(t,k;a)$ and $u(t,k;a_{0})$ to represent the solution corresponding to the Cauchy problems with coefficients $a(t)$ and $a_{0}$, respectively.
\subsection{Existence and uniqueness}\label{extnun}
Taking the $\mathcal{H}_{\hbar,V}$-Fourier transform of the Cauchy problem \eqref{mainpde} with respect to the variable $k\in\hbar\mathbb{Z}^{n}$ and using \eqref{symbol}, we obtain the collection  of Cauchy problems with Fourier coefficients: 
\begin{equation}\label{maintranspde}
	\left\{\begin{array}{l}
		\mathbb{D}_{(g)} \widehat{u}(t, \xi;a)+\lambda_{\xi}a(t) \widehat{u}(t, \xi;a)=\widehat{f}(t, \xi), \quad t\in(0,T],  \xi \in  \mathcal{I}_{\hbar}, \\
		\widehat{u}(0, \xi;a)=\widehat{u}_{0}(\xi), \quad \xi \in  \mathcal{I}_{\hbar},
	\end{array}\right.
\end{equation}
that is,
\begin{equation}\label{modtrans}
	\left\{\begin{array}{l}
		\mathbb{D}_{(g)} \widehat{u}(t, \xi;a)+\lambda_{\xi}a_{1}\widehat{u}(t, \xi;a)=\widehat{f}(t,\xi)+\lambda_{\xi}\left(a_{1}-a(t)\right)\widehat{u}(t, \xi;a),\quad t\in(0,T],  \\
		\widehat{u}(0, \xi;a)=\widehat{u}_{0}(\xi), \quad \xi \in  \mathcal{I}_{\hbar},
	\end{array}\right.
\end{equation}
where $a_{1}$ is a positive constant given by $$a_{1}:=\|a\|_{C([0,T])}=\sup\{|a(t)|:t\in[0,T]\}>0.$$ Since $\widehat{f}(t,\xi)+\lambda_{\xi}\left(a_{1}-a(t)\right)\widehat{u}(t, \xi;a)\in C([0,T])$, by using Remark \ref{3rdremark} for the Cauchy problem \eqref{modtrans}, we obtain a unique generalized solution  $\widehat{u}(t, \xi;a)\in C([0,T])$ which can be represented as
\begin{multline}\label{soluuxi}
	\widehat{u}(t, \xi;a)=\widehat{u}_{0}(\xi)\widehat{u}_{\lambda_{\xi}a_{1}}(t,\xi)-\\
	\frac{1}{\lambda_{\xi}a_{1}}\int^{t}_{0}\left[\widehat{f}(s,\xi)+\lambda_{\xi}\left(a_{1}-a(s)\right)\widehat{u}(s, \xi;a)\right]\frac{\partial}{\partial s}\widehat{u}_{\lambda_{\xi}a_{1}}(t-s,\xi)\mathrm{d}s,
\end{multline}
where the function $\widehat{u}_{\lambda_{\xi}a_{1}}(t,\xi)$ is the solution of the corresponding homogeneous equation of the Cauchy problem \eqref{modtrans} with initial Cauchy data $\widehat{u}(0, \xi;a)\equiv1$. Furthermore, the solution $u(t,k)=u(t,k;a)$ can be expressed in the form \eqref{specrep}, since the set of eigenfunctions $\{u_{\xi}\}_{\xi\in\mathcal{I}_{\hbar}}$ in \eqref{eigen} forms an orthonormal basis of $\ell^{2}(\hbar\mathbb{Z}^{n})$.
The uniqueness of solution $u(t,k)$ follows from the uniqueness of solution $\widehat{u}(t,\xi;a)$ in \eqref{soluuxi}. This completes the existence and uniqueness of solution $u(t,k)$ in Theorem \ref{mthm1}.
\subsection{Regularity estimates}
In order to establish the regularity estimate \eqref{reqesttt}, we split the Cauchy problem \eqref{mainpde} into homogeneous and the non-homogeneous Cauchy problems. Then, using the technique followed by the authors in \cite{2023niyaz2} and \cite{2023niyaz}, we will establish the regularity estimates for these Cauchy problems. Consequently, we will obtain the required estimates \eqref{reqesttt}.

Before beginning the proof, the following introductory lemma is required, which establishes the regularity of the solution $\widehat{u}(t,\xi; a)$.
\begin{lemma}\label{regularity}
	Let $\xi\in\mathcal{I}_{\hbar}$ and  $\widehat{u}(t,\xi ;a)$ be the unique solution of the Cauchy problem \eqref{maintranspde}. Then the following statements hold:
	\begin{enumerate}
		\item if $\widehat{f}(t,\xi) \leq 0$ on $[0, T]$ and $\widehat{u}_{0}(\xi) \leq 0$, then we have $\widehat{u}(t,\xi ; a) \leq 0$ on $[0, T]$,
		\item if $\widehat{f}(t,\xi) \geq 0$ on $[0, T]$ and $\widehat{u}_{0}(\xi) \geq 0$, then we have $\widehat{u}(t,\xi ; a) \geq 0$ on $[0, T]$.
	\end{enumerate}
\end{lemma}
\begin{proof} 
	The proof follows from Remark \ref{4thremark}.
\end{proof}
Let $u^{u_{0}}(t,k;a)$ and $u^{f}(t,k;a)$ be the solutions of the Cauchy problems
\begin{equation}\label{homeq1}
	\left\{\begin{array}{l}
		\mathbb{D}_{(g)} u^{u_{0}}(t, k;a)+a(t)\mathcal{H}_{\hbar,V} u^{u_{0}}(t, k;a)=0, \quad t \in(0, T],k\in \hbar \mathbb{Z}^n, \\
		u^{u_{0}}(0, k;a)=u_{0}(k), \quad k \in \hbar \mathbb{Z}^n,
	\end{array}\right.
\end{equation}
and
\begin{equation}\label{homeq2}
	\left\{\begin{array}{l}
		\mathbb{D}_{(g)} u^{f}(t, k;a)+a(t)\mathcal{H}_{\hbar,V} u^{f}(t, k;a)=f(t, k), \quad t\in(0, T],k\in \hbar \mathbb{Z}^n, \\
		u^{f}(0, k;a)=0, \quad k \in \hbar \mathbb{Z}^n,
	\end{array}\right.
\end{equation}
respectively. In the light of Theorem \ref{mthm1}, the solution $u(t,k;a)$ of the Cauchy problem \eqref{mainpde} can be written as
\begin{equation}\label{eqnnn}
	u(t,k;a)=u^{u_{0}}(t,k;a)+u^{f}(t,k;a),\quad t\in(0,T], k\in\hbar\mathbb{Z}^{n}.
\end{equation}
%Furthermore, using the spectral representation \eqref{specrep} the solution $\widehat{u}(t,\xi;a)$ of the Cauchy problem 
%\eqref{maintranspde} can be written as
%\begin{equation}
%	\widehat{u}(t,\xi;a)=\widehat{u}^{u_{0}}(t,\xi;a)+\widehat{u}^{f}(t,\xi;a),\quad t\in(0,T], k\in\hbar\mathbb{Z}^{n},
%\end{equation}
%where $\widehat{u}^{u_{0}}(t,\xi;a)$ and $\widehat{u}^{f}(t,\xi;a)$ are the  solutions of Cauchy problems
%\begin{equation}
%	\left\{\begin{array}{l}
	%		\mathbb{D}_{(g)} \widehat{u}^{u_{0}}(t, \xi;a)+\lambda_{\xi}a(t) \widehat{u}^{u_{0}}(t, \xi;a)=0, \quad t\in(0,T], \xi\in\mathcal{I}_{\hbar},\\
	%		\widehat{u}^{u_{0}}(0, \xi;a)=\widehat{u}_{0}(\xi), \quad \xi \in  \mathcal{I}_{\hbar},
	%	\end{array}\right.
%\end{equation}
%and
%\begin{equation}
%	\left\{\begin{array}{l}
	%		\mathbb{D}_{(g)} \widehat{u}^{f}(t, \xi;a)+\lambda_{\xi}a(t) \widehat{u}^{f}(t, \xi;a)=\widehat{f}(t, \xi), \quad t\in(0,T], \xi\in\mathcal{I}_{\hbar}, \\
	%		\widehat{u}^{f}(0, \xi;a)=0, \quad \xi \in  \mathcal{I}_{\hbar},
	%	\end{array}\right.
%\end{equation}
In the following Lemma \ref{unormm} and \ref{fnormm}, we establish the regularity estimate for the homogeneous Cauchy problem \eqref{homeq1} and non-homogeneous Cauchy problem \eqref{homeq2}, respectively.
\begin{lemma}\label{unormm}
	For $s\in\mathbb{R}$, let  $u_{0}\in \mathrm{H}^{2+s}_{\mathcal{H}_{\hbar,V}}$. If we assume that  $a\in C([0,T])$  is positive, then the solution $u^{u_{0}}(t,k;a)$ of the Cauchy problem \eqref{homeq1} satisfies the estimate 
	\begin{equation*}
		\|u^{u_{0}}(t,\cdot)\|_{\mathrm{H}^{2+s}_{\mathcal{H}_{\hbar,V}}}+	\|\mathbb{D}_{(g)}u^{u_{0}}(t,\cdot)\|_{\mathrm{H}^{s}_{\mathcal{H}_{\hbar,V}}}\leq \left(1+\|a\|_{C([0,T])}\right) \|u_{0}\|_{\mathrm{H}^{2+s}_{\mathcal{H}_{\hbar,V}}},
	\end{equation*}	
	for all $t\in[0,T]$.
\end{lemma}
\begin{proof} Consider the Cauchy problems
	\begin{equation}\label{maintranspde1}
		\mathbb{D}_{(g)} \widehat{u}^{u_{0}}(t, \xi;a)+\lambda_{\xi}a(t) \widehat{u}^{u_{0}}(t, \xi;a)=0, \quad
		\widehat{u}^{u_{0}}(0, \xi;a)=\widehat{u}_{0}(\xi), \quad \xi \in  \mathcal{I}_{\hbar},
	\end{equation}
	and 
	\begin{equation}\label{maintranspde2}
		\mathbb{D}_{(g)} \widehat{u}^{u_{0}}(t, \xi;a_{0})+\lambda_{\xi}a_{0} \widehat{u}^{u_{0}}(t, \xi;a_{0})=0, \quad
		\widehat{u}^{u_{0}}(0, \xi;a_{0})=\widehat{u}_{0}(\xi), \quad \xi \in  \mathcal{I}_{\hbar}.
	\end{equation}
	Applying Lemma \ref{regularity} to the Cauchy problem \eqref{maintranspde1}, we  deduce that  $\widehat{u}^{u_{0}}(t, \xi;a)$ holds the same sign as $\widehat{u}_{0}(\xi)$ and  applying the Lemma \ref{homthm} to the Cauchy problem \eqref{maintranspde2}, we get
	\begin{equation}\label{u0solution}
		\widehat{u}^{u_{0}}(t, \xi;a_{0})=
		\widehat{u}_{0}(\xi)\widehat{u}^{u_{0}}_{\lambda_{\xi}a_{0}}(t,\xi),\quad t\in(0,T],\xi\in\mathcal{I}_{\hbar},
	\end{equation}
	where $\widehat{u}^{u_{0}}_{\lambda_{\xi}a_{0}}(t,\xi)$ is the solution of the Cauchy problem \eqref{maintranspde2} with initial Cauchy data $\widehat{u}^{u_{0}}(0,\xi;a_{0})\equiv 1$.
	From the equations \eqref{maintranspde1} and  \eqref{maintranspde2}, denoting $w^{u_{0}}(t,\xi;a):=\widehat{u}^{u_{0}}(t, \xi;a_{0})-\widehat{u}^{u_{0}}(t, \xi;a)$, we get
	\begin{equation}\label{maintranspde3}
		\mathbb{D}_{(g)}w^{u_{0}}(t, \xi;a)+\lambda_{\xi}a_{0} w^{u_{0}}(t, \xi;a)=\lambda_{\xi}\widehat{u}^{u_{0}}(t, \xi;a)(a(t)-a_{0}), \quad
		w^{u_{0}}(0, \xi;a)=0, 
	\end{equation}
	for all $\xi\in\mathcal{I}_{\hbar}$. Now using the regularity of  $\widehat{u}^{u_{0}}(t, \xi;a)$ and  Lemma \ref{regularity} for the Cauchy problem \eqref{maintranspde3}, we deduce that the solution $w^{u_{0}}(t, \xi;a)$ also holds the same sign as $\widehat{u}_{0}(\xi)$, since $\lambda_{\xi}\left(a(t)-a_{0}\right)\geq0$ for all $t\in[0,T]$. Now combining the regularity of solutions $w^{u_{0}}(t, \xi;a)$ and $\widehat{u}^{u_{0}}(t, \xi;a_{0})$ with equation \eqref{u0solution}, we get
	\begin{equation*}
		\left\{\begin{array}{l}
			0 \leq \widehat{u}^{u_{0}}(t, \xi;a) \leq \widehat{u}^{u_{0}}(t, \xi;a_{0})=\widehat{u}_{0}(\xi)\widehat{u}^{u_{0}}_{\lambda_{\xi}a_{0}}(t,\xi),\quad \text { if } \widehat{u}_{0}(\xi) \geq 0, \\
			\widehat{u}_{0}(\xi)\widehat{u}^{u_{0}}_{\lambda_{\xi}a_{0}}(t,\xi)=\widehat{u}^{u_{0}}(t, \xi;a_{0}) \leq \widehat{u}^{u_{0}}(t, \xi;a) \leq 0,\quad \text { if }\widehat{u}_{0}(\xi) \leq 0.
		\end{array}\right.
	\end{equation*}
	Recalling Remark \ref{2ndremark}, i.e., $0\leq \widehat{u}^{u_{0}}_{\lambda_{\xi}a_{0}}(t,\xi)\leq 1$,  for all $(t,\xi)\in[0,T]\times\mathcal{I}_{\hbar}$, we obtain
	\begin{equation}\label{u0est}
		|\widehat{u}^{u_{0}}(t, \xi;a)\leq|\widehat{u}_{0}(\xi)\widehat{u}^{u_{0}}_{\lambda_{\xi}a_{0}}(t,\xi)|\leq |\widehat{u}_{0}(\xi)|,\quad \text{for all } t\in(0,T],\xi\in\mathcal{I}_{\hbar}.
	\end{equation}
	Additionally, using the above estimate and the Plancherel formula \eqref{fnorm}, we get 
	\begin{equation}\label{est11}
		\left\|u^{u_{0}}(t,\cdot ; a)\right\|^{2}_{\mathrm{H}^{2+s}_{\mathcal{H}_{\hbar,V}}}=\sum_{\xi \in \mathcal{I}_{\hbar}}\left(1+\lambda_{\xi}\right)^{2+s} \left|\widehat{u}^{u_{0}}(t,\xi ; a)\right|^{2}
		\leq  \|u_{0}\|^{2}_{\mathrm{H}^{2+s}_{\mathcal{H}_{\hbar,V}}}.
	\end{equation}
	Similarily, using the Plancherel formula \eqref{fnorm} and equation \eqref{homeq1} with estimate \eqref{u0est}, we get \begin{eqnarray}\label{est22}
		\left\|\mathbb{D}_{(g)}u^{u_{0}}(t,\cdot ; a)\right\|^{2}_{\mathrm{H}^{s}_{\mathcal{H}_{\hbar,V}}}&=&\|-a(t)\mathcal{H}_{\hbar,V}u^{u_{0}}(t,\cdot ; a)\|^{2}_{\mathrm{H}^{s}_{\mathcal{H}_{\hbar,V}}}\nonumber\\
		&\leq&\|a\|^{2}_{C([0,T])}\sum_{\xi \in \mathcal{I}_{\hbar}}\left(1+\lambda_{\xi}\right)^{s} \left|\lambda_{\xi}\widehat{u}^{u_{0}}(t,\xi ; a)\right|^{2}\nonumber\\
		&\leq& \|a\|^{2}_{C([0,T])}\sum_{\xi \in \mathcal{I}_{\hbar}}\left(1+\lambda_{\xi}\right)^{2+s} \left|\widehat{u}_{0}(\xi)\right|^{2} \nonumber\\
		&=& \|a\|^{2}_{C([0,T])}\|u_{0}\|^{2}_{\mathrm{H}^{2+s}_{\mathcal{H}_{\hbar,V}}},
	\end{eqnarray}
	for all $t\in[0,T].$ Putting together \eqref{est11} and \eqref{est22}, we get
	\begin{equation*}
		\|u^{u_{0}}(t,\cdot)\|_{\mathrm{H}^{2+s}_{\mathcal{H}_{\hbar,V}}}+	\|\mathbb{D}_{(g)}u^{u_{0}}(t,\cdot)\|_{\mathrm{H}^{s}_{\mathcal{H}_{\hbar,V}}}\leq \left(1+\|a\|_{C([0,T])}\right) \|u_{0}\|_{\mathrm{H}^{2+s}_{\mathcal{H}_{\hbar,V}}},
	\end{equation*}	
	for all $t\in[0,T]$.
	This completes the proof.
\end{proof}
\begin{lemma}\label{fnormm}
	For $s\in\mathbb{R}$, let  $f\in C([0,T];\mathrm{H}^{2+s}_{\mathcal{H}_{\hbar,V}})$. If we assume that  $a\in C([0,T])$ is positive, then the solution $u^{f}(t,k;a)$ of the Cauchy problem \eqref{homeq2} satisfies the estimate 
	\begin{multline*}
		\|u^{f}(t,\cdot)\|_{\mathrm{H}^{2+s}_{\mathcal{H}_{\hbar,V}}}+	\|\mathbb{D}_{(g)}u^{f}(t,\cdot)\|_{\mathrm{H}^{s}_{\mathcal{H}_{\hbar,V}}}\leq \\	(1+\|a\|_{C([0,T])})\max\left\{1,\frac{1}{\Lambda a_{0}}\right\} \|f\|_{C([0,T];\mathrm{H}^{2+s}_{\mathcal{H}_{\hbar,V}})},
	\end{multline*}	
	for all $t\in[0,T]$, where $\Lambda$ is given by \eqref{lambda}.
\end{lemma}
\begin{proof} Consider the Cauchy problems
	\begin{equation}\label{maintranspde11}
		\mathbb{D}_{(g)} \widehat{u}^{f}(t, \xi;a)+\lambda_{\xi}a(t) \widehat{u}^{f}(t, \xi;a)=\widehat{f}(t,\xi), \quad
		\widehat{u}^{f}(0, \xi;a)=0, \quad \xi \in  \mathcal{I}_{\hbar},
	\end{equation}
	and 
	\begin{equation}\label{maintranspde22}
		\mathbb{D}_{(g)} \widehat{u}^{f}(t, \xi;a_{0})+\lambda_{\xi}a_{0} \widehat{u}^{f}(t, \xi;a_{0})=\widehat{f}(t,\xi), \quad
		\widehat{u}^{f}(0, \xi;a_{0})=0, \quad \xi \in  \mathcal{I}_{\hbar}.
	\end{equation}
	If we write 
	\begin{equation*}
		\widehat{f}(t,\xi)=\mathrm{Re}     \widehat{f}(t,\xi)+i\mathrm{Im}    \widehat{f}(t,\xi),
	\end{equation*}
	then we can denote 
	\begin{equation*}
		\widehat{f}^{\max}_{\mathrm{Re}}(t,\xi):=\max\{0,\mathrm{Re}     \widehat{f}(t,\xi)\}=\frac{\mathrm{Re}     \widehat{f}(t,\xi)+|\mathrm{Re}     \widehat{f}(t,\xi)|}{2}\geq 0,     
	\end{equation*}
	\begin{equation*}
		\widehat{f}^{\min}_{\mathrm{Re}}(t,\xi):=\min\{0,\mathrm{Re}     \widehat{f}(t,\xi)\}=\frac{\mathrm{Re}     \widehat{f}(t,\xi)-|\mathrm{Re}     \widehat{f}(t,\xi)|}{2}\leq 0,     
	\end{equation*}
	\begin{equation*}
		\widehat{f}^{\max}_{\mathrm{Im}}(t,\xi):=\max\{0,\mathrm{Im}     \widehat{f}(t,\xi)\}=\frac{\mathrm{Im}     \widehat{f}(t,\xi)+|\mathrm{Im}     \widehat{f}(t,\xi)|}{2}\geq 0,     
	\end{equation*}
	and
	\begin{equation*}
		\widehat{f}^{\min}_{\mathrm{Im}}(t,\xi):=\min\{0,\mathrm{Im}     \widehat{f}(t,\xi)\}=\frac{\mathrm{Im}     \widehat{f}(t,\xi)-|\mathrm{Im}     \widehat{f}(t,\xi)|}{2}\leq 0.   
	\end{equation*}
	In view of the fact that $\mathrm{Re}\widehat{f}=\widehat{f}^{\max}_{\mathrm{Re}}+\widehat{f}^{\min}_{\mathrm{Re}}$ and $\mathrm{Im}\widehat{f}=\widehat{f}^{\max}_{\mathrm{Im}}+\widehat{f}^{\min}_{\mathrm{Im}}$, the Cauchy problems \eqref{maintranspde11} and \eqref{maintranspde22} can be split into the following Cauchy problems:
	\begin{equation}\label{fmaxa:Re}
		\mathbb{D}_{(g)} \widehat{u}^{f,\max}_{\mathrm{Re}}(t, \xi;a)+\lambda_{\xi}a(t) \widehat{u}^{f,\max}_{\mathrm{Re}}(t, \xi;a)=\widehat{f}^{\max}_{\mathrm{Re}}(t,\xi), \quad
		\widehat{u}^{f,\max}_{\mathrm{Re}}(0, \xi;a)=0,
	\end{equation}
	\begin{equation}\label{fmaxa0:Re}
		\mathbb{D}_{(g)} \widehat{u}^{f,\max}_{\mathrm{Re}}(t, \xi;a_{0})+\lambda_{\xi}a_{0} \widehat{u}^{f,\max}_{\mathrm{Re}}(t, \xi;a_{0})=\widehat{f}^{\max}_{\mathrm{Re}}(t,\xi), \quad
		\widehat{u}^{f,\max}_{\mathrm{Re}}(0, \xi;a_{0})=0,
	\end{equation}
	\begin{equation}\label{fmina:Re}
		\mathbb{D}_{(g)} \widehat{u}^{f,\min}_{\mathrm{Re}}(t, \xi;a)+\lambda_{\xi}a(t) \widehat{u}^{f,\min}_{\mathrm{Re}}(t, \xi;a)=\widehat{f}^{\min}_{\mathrm{Re}}(t,\xi), \quad
		\widehat{u}^{f,\min}_{\mathrm{Re}}(0, \xi;a)=0,
	\end{equation}
	\begin{equation}\label{fmina0:Re}		\mathbb{D}_{(g)} \widehat{u}^{f,\min}_{\mathrm{Re}}(t, \xi;a_{0})+\lambda_{\xi}a_{0} \widehat{u}^{f,\min}_{\mathrm{Re}}(t, \xi;a_{0})=\widehat{f}^{\min}_{\mathrm{Re}}(t,\xi), \quad
		\widehat{u}^{f,\min}_{\mathrm{Re}}(0, \xi;a_{0})=0,
	\end{equation}
	\begin{equation}\label{fmaxa:Im}
		\mathbb{D}_{(g)} \widehat{u}^{f,\max}_{\mathrm{Im}}(t, \xi;a)+\lambda_{\xi}a(t) \widehat{u}^{f,\max}_{\mathrm{Im}}(t, \xi;a)=\widehat{f}^{\max}_{\mathrm{Im}}(t,\xi), \quad
		\widehat{u}^{f,\max}_{\mathrm{Im}}(0, \xi;a)=0,
	\end{equation}
	\begin{equation}\label{fmaxa0:Im}
		\mathbb{D}_{(g)} \widehat{u}^{f,\max}_{\mathrm{Im}}(t, \xi;a_{0})+\lambda_{\xi}a_{0} \widehat{u}^{f,\max}_{\mathrm{Im}}(t, \xi;a_{0})=\widehat{f}^{\max}_{\mathrm{Im}}(t,\xi), \quad
		\widehat{u}^{f,\max}_{\mathrm{Im}}(0, \xi;a_{0})=0,
	\end{equation}
	\begin{equation}\label{fmina:Im}
		\mathbb{D}_{(g)} \widehat{u}^{f,\min}_{\mathrm{Im}}(t, \xi;a)+\lambda_{\xi}a(t) \widehat{u}^{f,\min}_{\mathrm{Im}}(t, \xi;a)=\widehat{f}^{\min}_{\mathrm{Im}}(t,\xi), \quad
		\widehat{u}^{f,\min}_{\mathrm{Im}}(0, \xi;a)=0,
	\end{equation}
	and
	\begin{equation}\label{fmina0:Im}		\mathbb{D}_{(g)} \widehat{u}^{f,\min}_{\mathrm{Im}}(t, \xi;a_{0})+\lambda_{\xi}a_{0} \widehat{u}^{f,\min}_{\mathrm{Im}}(t, \xi;a_{0})=\widehat{f}^{\min}_{\mathrm{Im}}(t,\xi), \quad
		\widehat{u}^{f,\min}_{\mathrm{Re}}(0, \xi;a_{0})=0.
	\end{equation}
	From the above equations, denoting 
	\begin{equation*}
		w^{f,\max}_{\mathrm{Re}}(t,\xi;a):=\widehat{u}^{f,\max}_{\mathrm{Re}}(t, \xi;a_{0})-\widehat{u}^{f,\max}_{\mathrm{Re}}(t, \xi;a),
	\end{equation*}
	\begin{equation*}
		w^{f,\min}_{\mathrm{Re}}(t,\xi;a):=\widehat{u}^{f,\min}_{\mathrm{Re}}(t, \xi;a_{0})-\widehat{u}^{f,\min}_{\mathrm{Re}}(t, \xi;a),
	\end{equation*}
	\begin{equation*}
		w^{f,\max}_{\mathrm{Im}}(t,\xi;a):=\widehat{u}^{f,\max}_{\mathrm{Im}}(t, \xi;a_{0})-\widehat{u}^{f,\max}_{\mathrm{Im}}(t, \xi;a),
	\end{equation*}
	and
	\begin{equation*}
		w^{f,\min}_{\mathrm{Im}}(t,\xi;a):=\widehat{u}^{f,\min}_{\mathrm{Im}}(t, \xi;a_{0})-\widehat{u}^{f,\min}_{\mathrm{Im}}(t, \xi;a),
	\end{equation*}
	we get
	\begin{equation}\label{wmax:Re}
		\left\{\begin{array}{l}
			\mathbb{D}_{(g)}w^{f,\max}_{\mathrm{Re}}(t, \xi;a)+\lambda_{\xi}a_{0} w^{f,\max}_{\mathrm{Re}}(t, \xi;a)=\lambda_{\xi}\widehat{u}^{f,\max}_{\mathrm{Re}}(t, \xi;a)(a(t)-a_{0}),\\
			w^{f,\max}_{\mathrm{Re}}(t, \xi;a)=0, \quad \xi \in  \mathcal{I}_{\hbar},t\in(0,T],
		\end{array}\right.
	\end{equation}
	\begin{equation}\label{wmin:Re}
		\left\{\begin{array}{l}
			\mathbb{D}_{(g)}w^{f,\min}_{\mathrm{Re}}(t, \xi;a)+\lambda_{\xi}a_{0} w^{f,\min}_{\mathrm{Re}}(t, \xi;a)=\lambda_{\xi}\widehat{u}^{f,\min}_{\mathrm{Re}}(t, \xi;a)(a(t)-a_{0}),\\
			w^{f,\min}_{\mathrm{Re}}(t, \xi;a)=0, \quad \xi \in  \mathcal{I}_{\hbar},t\in(0,T],
		\end{array}\right.
	\end{equation}
	\begin{equation}\label{wmax:Im}
		\left\{\begin{array}{l}
			\mathbb{D}_{(g)}w^{f,\max}_{\mathrm{Im}}(t, \xi;a)+\lambda_{\xi}a_{0} w^{f,\max}_{\mathrm{Im}}(t, \xi;a)=\lambda_{\xi}\widehat{u}^{f,\max}_{\mathrm{Im}}(t, \xi;a)(a(t)-a_{0}), \\
			w^{f,\max}_{\mathrm{Im}}(t, \xi;a)=0, \quad \xi \in  \mathcal{I}_{\hbar},t\in(0,T],
		\end{array}\right.
	\end{equation}
	and
	\begin{equation}\label{wmin:Im}
		\left\{\begin{array}{l}
			\mathbb{D}_{(g)}w^{f,\min}_{\mathrm{Im}}(t, \xi;a)+\lambda_{\xi}a_{0} w^{f,\min}_{\mathrm{Im}}(t, \xi;a)=\lambda_{\xi}\widehat{u}^{f,\min}_{\mathrm{Im}}(t, \xi;a)(a(t)-a_{0}),  \\
			w^{f,\min}_{\mathrm{Im}}(t, \xi;a)=0, \quad \xi \in  \mathcal{I}_{\hbar},t\in(0,T],.
		\end{array}\right.
	\end{equation}
	Now applying Lemma \ref{regularity} to the Cauchy problems \eqref{fmaxa:Re}, \eqref{fmina:Re}, \eqref{fmaxa:Im},   and \eqref{fmina:Im}, we deduce that
	\begin{equation}\label{cond1w}
		\left\{\begin{array}{l}
			\widehat{u}^{f, \max }_{\mathrm{Re}}(t,\xi ; a)\geq 0, \quad \text{since}~~ \widehat{f}^{\max}_{\mathrm{Re}}\geq 0;\\
			\widehat{u}^{f, \min }_{\mathrm{Re}}(t,\xi ; a)\leq 0, \quad \text{since}~~ \widehat{f}^{\min}_{\mathrm{Re}}\leq 0;\\
			\widehat{u}^{f, \max }_{\mathrm{Im}}(t,\xi ; a)\geq0, \quad \text{since}~~\widehat{f}^{\max}_{\mathrm{Im}}\geq 0; ~\text{and}\\
			\widehat{u}^{f, \min}_{\mathrm{Im}}(t,\xi ; a)\leq 0 \quad \text{since}~~ \widehat{f}^{\min}_{\mathrm{Im}}\leq 0;
		\end{array}\right.
	\end{equation}
	for all $t\in[0,T]$.
	Consequently, using the regularity of solutions in \eqref{cond1w}  with Lemma \ref{regularity} to the Cauchy problems \eqref{wmax:Re}-\eqref{wmin:Im}, we conclude that
	\begin{equation*}
		w^{f,\max}_{\mathrm{Re}}(t, \xi;a)\geq 0,\quad    	w^{f,\min}_{\mathrm{Re}}(t, \xi;a)\leq 0 ,
	\end{equation*}
	and \begin{equation}
		w^{f,\max}_{\mathrm{Im}}(t, \xi;a)\geq 0,\quad w^{f,\min}_{\mathrm{Im}}(t, \xi;a)\leq 0,
	\end{equation}
	for all $t\in[0,T]$, since $\lambda_{\xi}(a(t)-a_{0})\geq 0.$ This gives
	\begin{equation}\label{cond1}
		\left\{\begin{array}{l}
			0 \leq \widehat{u}^{f, \max }_{\mathrm{Re}}(t,\xi ; a) \leq \widehat{u}^{f, \max }_{\mathrm{Re}}\left(t ; a_{0}\right),\quad \text { on }[0, T], \\
			\widehat{u}^{f, \min }_{\mathrm{Re}}\left(t,\xi ; a_{0}\right) \leq \widehat{u}^{f, \min }_{\mathrm{Re}}(t ,\xi; a) \leq 0,\quad \text { on }[0, T],\\
			0 \leq \widehat{u}^{f, \max }_{\mathrm{Im}}(t,\xi ; a) \leq \widehat{u}^{f, \max }_{\mathrm{Im}}\left(t ; a_{0}\right),\quad \text { on }[0, T], \\
			\widehat{u}^{f, \min }_{\mathrm{Im}}\left(t,\xi ; a_{0}\right) \leq \widehat{u}^{f, \min }_{\mathrm{Im}}(t ,\xi; a) \leq 0,\quad \text { on }[0, T],
		\end{array}\right.
	\end{equation}
	for all $\xi\in\mathcal{I}_{\hbar}$.
	
	Now applying Lemma \ref{nonthm} to the Cauchy problems \eqref{maintranspde22}, \eqref{fmaxa0:Re}, \eqref{fmina0:Re}, \eqref{fmaxa0:Im}, and \eqref{fmina0:Im}, we get
	\begin{equation*}
		\widehat{u}^{f}(t,\xi;a_{0})=-\frac{1}{\lambda_{\xi}a_{0}} \int_{0}^{t} \widehat{f}(s,\xi)\frac{\partial}{\partial s}w_{\lambda_{\xi}a_{0}}(t-s,\xi)  \mathrm{d} s,
	\end{equation*}
	\begin{equation*}
		\widehat{u}^{f,\max}_{\mathrm{Re}}(t,\xi;a_{0})=-\frac{1}{\lambda_{\xi}a_{0}} \int_0^t \widehat{f}^{\max}_{\mathrm{Re}}(s,\xi)\frac{\partial}{\partial s}w_{\lambda_{\xi}a_{0}}(t-s,\xi)  \mathrm{d} s,
	\end{equation*}
	\begin{equation*}
		\widehat{u}^{f,\min}_{\mathrm{Re}}(t,\xi;a_{0})=-\frac{1}{\lambda_{\xi}a_{0}} \int_0^t \widehat{f}^{\min}_{\mathrm{Re}}(s,\xi)\frac{\partial}{\partial s}w_{\lambda_{\xi}a_{0}}(t-s,\xi)  \mathrm{d} s,
	\end{equation*}
	\begin{equation*}
		\widehat{u}^{f,\max}_{\mathrm{Im}}(t,\xi;a_{0})=-\frac{1}{\lambda_{\xi}a_{0}} \int_0^t \widehat{f}^{\max}_{\mathrm{Im}}(s,\xi)\frac{\partial}{\partial s}w_{\lambda_{\xi}a_{0}}(t-s,\xi)  \mathrm{d} s,
	\end{equation*}
	and
	\begin{equation*}
		\widehat{u}^{f,\min}_{\mathrm{Im}}(t,\xi;a_{0})=-\frac{1}{\lambda_{\xi}a_{0}} \int_0^t \widehat{f}^{\min}_{\mathrm{Im}}(s,\xi)\frac{\partial}{\partial s}w_{\lambda_{\xi}a_{0}}(t-s,\xi)  \mathrm{d} s,
	\end{equation*}
	respectively.
	
	Combining the facts that
	\begin{equation*}
		\widehat{u}^{f}(t,\xi;a)=\widehat{u}^{f,\max}_{\mathrm{Re}}(t,\xi;a)+\widehat{u}^{f,\min}_{\mathrm{Re}}(t,\xi;a)+i\left(\widehat{u}^{f,\max}_{\mathrm{Im}}(t,\xi;a)+\widehat{u}^{f,\min}_{\mathrm{Im}}(t,\xi;a)\right),
	\end{equation*}
	with the above equations and inequality \eqref{cond1}, we get
	\begin{eqnarray}\label{festi}
		&&|\widehat{u}^{f}(t,\xi;a)|^{2}\nonumber\\&\leq&|\widehat{u}^{f,\max}_{\mathrm{Re}}(t,\xi;a)+\widehat{u}^{f,\min}_{\mathrm{Re}}(t,\xi;a)|^{2}+|\widehat{u}^{f,\max}_{\mathrm{Im}}(t,\xi;a)+\widehat{u}^{f,\min}_{\mathrm{Im}}(t,\xi;a)|^{2}\nonumber\\
		&\leq&\left[\frac{1}{\lambda_{\xi}a_{0}} \int_0^t \left( |\widehat{f}^{\max}_{\mathrm{Re}}(s,\xi)|+|\widehat{f}^{\min}_{\mathrm{Re}}(s,\xi)|\right)\left|\frac{\partial}{\partial s}w_{\lambda_{\xi}a_{0}}(t-s,\xi) \right| \mathrm{d} s\right]^{2}\nonumber\\
		&+&\left[\frac{1}{\lambda_{\xi}a_{0}} \int_0^t \left( |\widehat{f}^{\max}_{\mathrm{Im}}(s,\xi)|+|\widehat{f}^{\min}_{\mathrm{Im}}(s,\xi)|\right)\left|\frac{\partial}{\partial s}w_{\lambda_{\xi}a_{0}}(t-s,\xi) \right| \mathrm{d} s\right]^{2}\nonumber\\
		&=&\left[\frac{1}{\lambda_{\xi}a_{0}} \int_0^t \left( \widehat{f}^{\max}_{\mathrm{Re}}(s,\xi)-\widehat{f}^{\min}_{\mathrm{Re}}(s,\xi)\right)\left|\frac{\partial}{\partial s}w_{\lambda_{\xi}a_{0}}(t-s,\xi) \right| \mathrm{d} s\right]^{2}~~ (\because \widehat{f}^{\min}_{\mathrm{Re}}\leq 0)\nonumber\\
		&+&\left[\frac{1}{\lambda_{\xi}a_{0}} \int_0^t \left( \widehat{f}^{\max}_{\mathrm{Im}}(s,\xi)-\widehat{f}^{\min}_{\mathrm{Im}}(s,\xi)\right)\left|\frac{\partial}{\partial s}w_{\lambda_{\xi}a_{0}}(t-s,\xi) \right| \mathrm{d} s\right]^{2}~ (\because \widehat{f}^{\min}_{\mathrm{Im}}\leq 0)\nonumber\\
		&=&\left[\frac{1}{\lambda_{\xi}a_{0}} \int_0^t  \left|\mathrm{Re}\widehat{f}(s,\xi)\right|\left|\frac{\partial}{\partial s}w_{\lambda_{\xi}a_{0}}(t-s,\xi) \right| \mathrm{d} s\right]^{2}\nonumber\\
		&+&\left[\frac{1}{\lambda_{\xi}a_{0}} \int_0^t \left|\mathrm{Im}\widehat{f}(s,\xi)\right|\left|\frac{\partial}{\partial s}w_{\lambda_{\xi}a_{0}}(t-s,\xi) \right| \mathrm{d} s\right]^{2}\nonumber\\
		&\leq&\left[\frac{1}{\lambda_{\xi}a_{0}} \sup_{s\in[0,T]}\left\{\left|\mathrm{Re}\widehat{f}(s,\xi)\right|\right\}\int_0^t  \left|\frac{\partial}{\partial s}w_{\lambda_{\xi}a_{0}}(t-s,\xi) \right| \mathrm{d} s\right]^{2}\nonumber\\
		&+&\left[\frac{1}{\lambda_{\xi}a_{0}} \sup_{s\in[0,T]}\left\{\left|\mathrm{Im}\widehat{f}(s,\xi)\right|\right\}\int_0^t \left|\frac{\partial}{\partial s}w_{\lambda_{\xi}a_{0}}(t-s,\xi) \right| \mathrm{d} s\right]^{2}\nonumber\\
		&=&\left[\frac{1}{\lambda_{\xi}a_{0}} \sup_{s\in[0,T]}\left\{\left|\mathrm{Re}\widehat{f}(s,\xi)\right|\right\}\left(1-w_{\lambda_{\xi}a_{0}}(t,\xi)\right)\right]^{2}\nonumber\\
		&+&\left[\frac{1}{\lambda_{\xi}a_{0}} \sup_{s\in[0,T]}\left\{\left|\mathrm{Im}\widehat{f}(s,\xi)\right|\right\}\left(1-w_{\lambda_{\xi}a_{0}}(t,\xi)\right)\right]^{2}\quad (\because w_{\lambda_{\xi}a_{0}}(0,\xi)=1)\nonumber\\
		&\leq&\frac{1}{(\lambda_{\xi}a_{0})^{2}}\sup_{s\in[0,T]}\left\{\left|\mathrm{Re}\widehat{f}(s,\xi)\right|^{2}+\left|\mathrm{Im}\widehat{f}(s,\xi)\right|^{2}\right\}\nonumber\\
		&=&\frac{1}{(\lambda_{\xi}a_{0})^{2}}\sup_{s\in[0,T]}\left\{\left|\widehat{f}(s,\xi)\right|^{2}\right\}\nonumber\\
		&\leq& \frac{\|\widehat{f}(\cdot,\xi)\|^{2}_{C([0,T])}}{(\lambda_{\xi}a_{0})^2},
	\end{eqnarray}
	for all $t\in[0,T]$ and $\xi\in\mathcal{I}_{\hbar}$, where the last inequality follows from the fact that $0\leq w_{\lambda_{\xi}a_{0}}(t,\xi)\leq 1$, for all $(t,\xi)\in[0,T]\times\mathcal{I}_{\hbar}$, see Remark \ref{2ndremark}.

	Further, using the Plancherel formula \eqref{fnorm} with the estimate \eqref{festi}, we get
	\begin{eqnarray}\label{est1}
		\left\|u^{f}(t,\cdot ; a)\right\|^{2}_{\mathrm{H}^{2+s}_{\mathcal{H}_{\hbar,V}}}&=&\sum_{\xi \in \mathcal{I}_{\hbar}}\left(1+\lambda_{\xi}\right)^{2+s} \left|\widehat{u}^{f}(t,\xi ; a)\right|^{2}\nonumber\\
		&\leq& \sum_{\xi \in \mathcal{I}_{\hbar}}\left(1+\lambda_{\xi}\right)^{2+s} \frac{\|\widehat{f}(\cdot,\xi)\|^{2}_{C([0,T])}}{\left(\lambda_{\xi}a_{0}\right)^{2}}\nonumber\\
		&\leq&\frac{1}{(\Lambda a_{0})^{2}} \|f\|^{2}_{C([0,T];\mathrm{H}^{2+s}_{\mathcal{H}_{\hbar,V}})}, 
	\end{eqnarray}
	where $\Lambda$ is given by \eqref{lambda}. Similarly, combining the Plancherel formula \eqref{fnorm} and equation \eqref{homeq2} with the estimate \eqref{festi}, we get
	\begin{eqnarray}\label{est2}
		&&\left\|\mathbb{D}_{(g)}u^{f}(t,\cdot ; a)\right\|^{2}_{\mathrm{H}^{s}_{\mathcal{H}_{\hbar,V}}}\nonumber\\&=&\|-a(t)\mathcal{H}_{\hbar,V}u^{f}(t,\cdot ; a)+f(t,\cdot)\|^{2}_{\mathrm{H}^{s}_{\mathcal{H}_{\hbar,V}}}\nonumber\\
		&\leq&\|a\|^{2}_{C([0,T])}\sum_{\xi \in \mathcal{I}_{\hbar}}\left(1+\lambda_{\xi}\right)^{s} \left|\lambda_{\xi}\widehat{u}^{f}(t,\xi ; a)\right|^{2}+\|f(t,\cdot)\|^{2}_{\mathrm{H}^{s}_{\mathcal{H}_{\hbar,V}}}\nonumber\\
		&\leq& \|a\|^{2}_{C([0,T])}\sum_{\xi \in \mathcal{I}_{\hbar}}\left(1+\lambda_{\xi}\right)^{2+s} \frac{\|\widehat{f}(\cdot,\xi)\|^{2}_{C([0,T])}}{\left(\lambda_{\xi}a_{0}\right)^{2}}+\|f\|^{2}_{C([0,T];\mathrm{H}^{2+s}_{\mathcal{H}_{\hbar,V}})} \nonumber\\
		&\leq& \left[\frac{\|a\|^{2}_{C([0,T])}}{(\Lambda a_{0})^{2}}+1\right]\|f\|^{2}_{C([0,T];\mathrm{H}^{2+s}_{\mathcal{H}_{\hbar,V}})}\nonumber\\&=&\left[(1+\|a\|_{C([0,T])})\max\left\{1,\frac{1}{\Lambda a_{0}}\right\}\right]^{2}\|f\|^{2}_{C([0,T];\mathrm{H}^{2+s}_{\mathcal{H}_{\hbar,V}})}, 
	\end{eqnarray}
	for all $t\in[0,T]$. Putting together estimates \eqref{est1} and \eqref{est2},
	\begin{multline*}
		\|u^{f}(t,\cdot)\|_{\mathrm{H}^{2+s}_{\mathcal{H}_{\hbar,V}}}+	\|\mathbb{D}_{(g)}u^{f}(t,\cdot)\|_{\mathrm{H}^{s}_{\mathcal{H}_{\hbar,V}}}\leq\\ 
		(1+\|a\|_{C([0,T])})\max\left\{1,\frac{1}{\Lambda a_{0}}\right\} \|f\|_{C([0,T];\mathrm{H}^{2+s}_{\mathcal{H}_{\hbar,V}})},
	\end{multline*}
	for all $t\in[0,T]$. This completes the proof.
\end{proof}
Now, we are in a position to prove the well-posedness theorem for the Cauchy problem (\ref{mainpde}) with the regular coefficient and the source term. 
\begin{proof}[Proof of Theorem \ref{mthm1}] Combining the equality \eqref{eqnnn} with Lemmas \ref{unormm} and \ref{fnormm} we obtain the following regularity estimate
	\begin{multline*}
		\|u(t,\cdot)\|_{\mathrm{H}^{2+s}_{\mathcal{H}_{\hbar,V}}}+	\|\mathbb{D}_{(g)}u(t,\cdot)\|_{\mathrm{H}^{s}_{\mathcal{H}_{\hbar,V}}}\leq\\ (1+\|a\|_{C([0,T])})\max\left\{1,\frac{1}{\Lambda a_{0}}\right\}\left[\|u_{0}\|_{\mathrm{H}^{2+s}_{\mathcal{H}_{\hbar,V}}}+\|f\|_{C([0,T];\mathrm{H}^{2+s}_{\mathcal{H}_{\hbar,V}})}\right].    
	\end{multline*}
	This completes the proof.
\end{proof}
In the above result, we proved the well-posedness to (\ref{mainpde}) with regular coefficient in the Sobolev-type spaces. Now, we prove the existence and uniqueness of very weak solutions in the case of distributional coefficients. 
\begin{proof}[Proof of Theorem \ref{ext}]
	\sloppy	Let $a_{\varepsilon}$ be the $C$-moderate and $f_{\varepsilon}$ be the $C([0, T] ; \mathrm{H}_{\mathcal{H}_{\hbar,V}}^{2+s})$-moderate regularization of the coefficient $a$ and the source term $f$, respectively. Now, fix $\varepsilon\in(0,1]$, and consider the regularized problem 
	\begin{equation}\label{wkcauchy}
		\left\{\begin{array}{l}
			\mathbb{D}_{(g)} u_{\varepsilon}(t, k)+a_{\varepsilon}(t)\mathcal{H}_{\hbar,V}u_{\varepsilon}(t, k)=f_{\varepsilon}(t, k),\quad(t, k) \in(0, T] \times \hbar\mathbb{Z}^{n}, \\
			u_{\varepsilon}(0, k)=u_{0}(k),\quad k \in \hbar\mathbb{Z}^{n}.
		\end{array}\right.
	\end{equation}
	Since $a\in \mathcal{D}^{\prime}(\mathbb{R})$, by the structure theorem for compactly supported distributions, there exists $L_{1}\in\mathbb{N}$ such that
	\begin{equation}\label{aest}
		|a_{\varepsilon}(t)|\lesssim \omega(\varepsilon)^{-L_{1}},\quad \text{for all }t\in[0,T],
	\end{equation}
	where $\omega(\varepsilon)\to 0$ as $\varepsilon\to 0$. Since $a\geq a_{0}>0$, it follows that
	\begin{equation}
		a_{\varepsilon}(t)=(a*\psi_{\omega(\varepsilon)})(t)=\langle a,\tau_{t}\tilde{\psi}_{\omega(\varepsilon)}\rangle\geq \tilde{a}_{0}>0,
	\end{equation}
	where $\tilde{\psi}(x)=\psi(-x), x \in \mathbb{R}$ and $\tau_t$ is a translation operator. Using Theorem \ref{mthm1} to the Cauchy problem \eqref{wkcauchy}, we deduce that
	\begin{multline}\label{reqestttt}
		\|u_{\varepsilon}(t,\cdot)\|_{\mathrm{H}^{2+s}_{\mathcal{H}_{\hbar,V}}}+	\|\mathbb{D}_{(g)}u_{\varepsilon}(t,\cdot)\|_{\mathrm{H}^{s}_{\mathcal{H}_{\hbar,V}}}\lesssim\\\left(1+\|a_{\varepsilon}\|_{C([0,T])}\right)\left[\|u_{0}\|_{\mathrm{H}^{2+s}_{\mathcal{H}_{\hbar,V}}}+\|f_{\varepsilon}\|_{C([0,T];\mathrm{H}^{2+s}_{\mathcal{H}_{\hbar,V}})}\right],       
	\end{multline}
	for all $t\in[0,T]$. Using the estimate \eqref{aest} and the inequality $1+x\leq e^x,x\in \mathbb{R}$, we get 
	\begin{equation}
		\|u_{\varepsilon}(t,\cdot)\|_{\mathrm{H}^{2+s}_{\mathcal{H}_{\hbar,V}}}+	\|\mathbb{D}_{(g)}u_{\varepsilon}(t,\cdot)\|_{\mathrm{H}^{s}_{\mathcal{H}_{\hbar,V}}}\lesssim e^{\omega(\varepsilon)^{-L_{1}}}\left[\|u_{0}\|_{\mathrm{H}^{2+s}_{\mathcal{H}_{\hbar,V}}}+\|f_{\varepsilon}\|_{C([0,T];\mathrm{H}^{2+s}_{\mathcal{H}_{\hbar,V}})}\right],     
	\end{equation}
	for all $t\in[0,T]$. Putting $\omega(\varepsilon)^{-L_{1}}=\log(\varepsilon^{-1})$, we obtain
	\begin{equation}\label{reqestttt0}
		\|u_{\varepsilon}(t,\cdot)\|_{\mathrm{H}^{2+s}_{\mathcal{H}_{\hbar,V}}}+	\|\mathbb{D}_{(g)}u_{\varepsilon}(t,\cdot)\|_{\mathrm{H}^{s}_{\mathcal{H}_{\hbar,V}}}\lesssim \varepsilon^{-1}\left[\|u_{0}\|_{\mathrm{H}^{2+s}_{\mathcal{H}_{\hbar,V}}}+\|f_{\varepsilon}\|_{C([0,T];\mathrm{H}^{2+s}_{\mathcal{H}_{\hbar,V}})}\right],   
	\end{equation}
	for all $t\in[0,T]$. Since $f_{\varepsilon}$ is $C([0,T];\mathrm{H}_{\mathcal{H}_{\hbar,V}}^{2+s})$-moderate, there exists $L_{2}\in \mathbb{N}$ such that
	\begin{equation}\label{afnorm}
		\|f_{\varepsilon}\|_{C([0,T];\mathrm{H}_{\mathcal{H}_{\hbar,V}}^{2+s})}\lesssim \varepsilon^{-L_{2}}.
	\end{equation}
	Combining the estimates \eqref{reqestttt0} and \eqref{afnorm}, we get
	\begin{equation*}
		\|u_{\varepsilon}\|_{C([0,T];\mathrm{H}_{\mathcal{H}_{\hbar,V}}^{2+s})}\lesssim \varepsilon^{-1-L_{2}}. 
	\end{equation*}
	Thus, we conclude that $u_{\varepsilon}$ is $C([0,T];\mathrm{H}_{\mathcal{H}_{\hbar,V}}^{2+s})$-moderate. This proves the existence of a very weak solution for the Cauchy problem \eqref{mainpde}.
	
	Now we will prove the uniqueness of a very weak solution. 
	Let $(u_{\varepsilon})_{\varepsilon}$	 and $(\tilde{u}_{\varepsilon})_{\varepsilon}$	 solve the Cauchy problems \eqref{reg1} and \eqref{reg2}, respectively. 
	Denoting $w_{\varepsilon}:=u_{\varepsilon}-\tilde{u}_{\varepsilon}$, we obtain
	\begin{equation}\label{weqn}
		\left\{\begin{array}{l}
			\mathbb{D}_{(g)} w_{\varepsilon}(t, k)+a_{\varepsilon}(t)\mathcal{H}_{\hbar,V}w_{\varepsilon}(t, k)=g_{\varepsilon}(t, k),\quad (t,k)\in (0,T]\times\hbar\mathbb{Z}^{n}, \\
			w_{\varepsilon}(0, k)=0,\quad k \in \hbar\mathbb{Z}^{n}, 
		\end{array}\right.
	\end{equation}
	where the source term $g_{\varepsilon}$ is given by
	\begin{equation}
		g_{\varepsilon}(t,k):=(\tilde{a}_{\varepsilon}-a_{\varepsilon})(t) \mathcal{H}_{\hbar,V}\tilde{u}_{\varepsilon}(t, k)+(f_{\varepsilon}-\tilde{f}_{\varepsilon})(t,k).
	\end{equation}
	Using the fact that  $(\tilde{a}_{\varepsilon}-a_{\varepsilon})_{\varepsilon}$   and $(f_{\varepsilon}-\tilde{f}_{\varepsilon})_{\varepsilon}$ are $C$-negligible and $C([0,T];\mathrm{H}_{\mathcal{H}_{\hbar,V}}^{2+s})$-negligible, respectively, it is evident that $g_{\varepsilon}$ is $C([0,T];\mathrm{H}_{\mathcal{H}_{\hbar,V}}^{2+s})$-negligible.  Since $w_{\varepsilon}(0,k)\equiv 0$,  it follows from Theorem \ref{mthm1}, estimate \eqref{aest} and the inequality $1+x\leq e^{x}$ that 
	\begin{equation*}
		\|w_{\varepsilon}(t,\cdot)\|_{\mathrm{H}^{2+s}_{\mathcal{H}_{\hbar,V}}}+	\|\mathbb{D}_{(g)}w_{\varepsilon}(t,\cdot)\|_{\mathrm{H}^{s}_{\mathcal{H}_{\hbar,V}}}\lesssim e^{\omega(\varepsilon)^{-L_{1}}}\|g_{\varepsilon}\|_{C([0,T];\mathrm{H}^{2+s}_{\mathcal{H}_{\hbar,V}})}.
	\end{equation*}
	Putting $\omega(\varepsilon)^{-L_{1}}=\log(\varepsilon^{-1})$ and using the fact that $g_{\varepsilon}$ is $C([0,T];\mathrm{H}^{2+s}_{\mathcal{H}_{\hbar,V}})$-negligible, it follows that
	\begin{equation*}
		\|w_{\varepsilon}(t,\cdot)\|_{\mathrm{H}^{2+s}_{\mathcal{H}_{\hbar,V}}}+	\|\mathbb{D}_{(g)}w_{\varepsilon}(t,\cdot)\|_{\mathrm{H}^{s}_{\mathcal{H}_{\hbar,V}}}\lesssim \varepsilon^{-1}\varepsilon^{-1+q}=\varepsilon^{q},\quad q\in \mathbb{N}.
	\end{equation*}
	Thus $w_{\varepsilon}$ is $C([0,T];\mathrm{H}^{2+s}_{\mathcal{H}_{\hbar,V}})$-negligible. This proves the uniqueness of a very weak solution and, hence, completes the proof.
\end{proof}
In the above result, we discussed the uniqueness of the very weak solution in the sense of Definition \ref{uniquedef}. Now, we prove the consistency of the very weak solution with the classical solution.
\begin{proof}[Proof of Theorem \ref{cnst}]
	Let $u$ be the classical solution of the Cauchy problem:
	\begin{equation}\label{cnst1}
		\left\{\begin{array}{l}
			\mathbb{D}_{(g)} u(t, k)+a(t) \mathcal{H}_{\hbar,V}u(t, k)=f(t,k),\quad(t, k) \in(0, T] \times \hbar\mathbb{Z}^{n}, \\
			u(0, k)=u_{0}(k),\quad k \in \hbar\mathbb{Z}^{n},
		\end{array}\right.
	\end{equation} 	
	and $(u_{\varepsilon})_{\varepsilon}$ be the very weak solution, i.e., $(u_{\varepsilon})_{\varepsilon}$ satisfies the following regularized Cauchy problem
	\begin{equation}\label{cnst2}
		\left\{\begin{array}{l}
			\mathbb{D}_{(g)} u_{\varepsilon}(t, k)+a_{\varepsilon}(t) \mathcal{H}_{\hbar,V}u_{\varepsilon}(t, k)=f_{\varepsilon}(t,k),\quad(t, k) \in(0, T] \times \hbar\mathbb{Z}^{n}, \\
			u_{\varepsilon}(0, k)=u_{0}(k),\quad k \in \hbar\mathbb{Z}^{n}.
		\end{array}\right.
	\end{equation} 
	In the light of \cite[Theorem 1.5.13]{Ruz:book}, we note that 
	\begin{equation}
		\left(a_{\varepsilon}-a\right)_{\varepsilon}\to 0  \text{ in } C([0,T]), \quad\varepsilon\to 0,  
	\end{equation}
	and 
	\begin{equation}
		\left(f-f_{\varepsilon}\right)_{\varepsilon}\to 0 \text{ in } C([0,T];\mathrm{H}^{2+s}_{\mathcal{H}_{\hbar,V}}),\quad \varepsilon\to 0,
	\end{equation}
	for $a\in C([0,T])$ and $f\in C([0,T];\mathrm{H}^{2+s}_{\mathcal{H}_{\hbar,V}})$, respectively. Further, we can rewrite the Cauchy problem \eqref{cnst1} as follows:
	\begin{equation}\label{cnst3}
		\left\{\begin{array}{l}
			\mathbb{D}_{(g)} \tilde{u}(t, k)+a_{\varepsilon}(t) \mathcal{H}_{\hbar,V}\tilde{u}(t, k)=f_{\varepsilon}(t,k)+g_{\varepsilon}(t,k),~~ (t, k) \in(0, T] \times \hbar\mathbb{Z}^{n}, \\
			\tilde{u}(0, k)=u_{0}(k),\quad k \in \hbar\mathbb{Z}^{n},
		\end{array}\right.
	\end{equation} 	
	where $g_{\varepsilon}$ is given by $$g_{\varepsilon}(t,k):=\left(a_{\varepsilon}-a\right)(t) \mathcal{H}_{\hbar,V} u(t, k)+\left(f-f_{\varepsilon}\right)(t,k),$$
	$g_{\varepsilon} \in C([0, T] ; \mathrm{H}_{\mathcal{H}_{\hbar,V}}^{2+s})$ and  $g_{\varepsilon}\to0$ in $C([0, T] ; \mathrm{H}_{\mathcal{H}_{\hbar,V}}^{2+s})$ as $\varepsilon \rightarrow 0$. Using the Cauchy problems, \eqref{cnst2} and \eqref{cnst3}, denoting  $w_{\varepsilon}:=u-u_{\varepsilon}$,  we obtain
	\begin{equation*}\label{cnst4}
		\left\{\begin{array}{l}
			\mathbb{D}_{(g)} w_{\varepsilon}(t, k)+a_{\varepsilon}(t) \mathcal{H}_{\hbar,V}w_{\varepsilon}(t, k)=g_{\varepsilon}(t,k),\quad(t, k) \in(0, T] \times \hbar\mathbb{Z}^{n}, \\
			w_{\varepsilon}(0, k)=0,\quad k \in \hbar\mathbb{Z}^{n}.
		\end{array}\right.
	\end{equation*}
	Since the coefficients and source term are regular enough, it follows from Theorem \ref{mthm1} that
	\begin{equation*}
		\|w_{\varepsilon}(t,\cdot)\|_{\mathrm{H}_{\mathcal{H}_{\hbar,V}}^{2+s}}+\|\mathbb{D}_{(g)}w_{\varepsilon}(t,\cdot)\|_{\mathrm{H}_{\mathcal{H}_{\hbar,V}}^{s}}\lesssim 	\|g_{\varepsilon}\|_{C([0,T];\mathrm{H}_{\mathcal{H}_{\hbar,V}}^{2+s})},
	\end{equation*}
	for all $t\in [0,T]$.
	Taking sup-norm over $[0,T]$ in the above estimate, we obtain
	\begin{equation*}
		\|w_{\varepsilon}\|_{C([0,T];\mathrm{H}_{\mathcal{H}_{\hbar,V}}^{2+s})}+\|\mathbb{D}_{(g)}w_{\varepsilon}\|_{C([0,T];\mathrm{H}_{\mathcal{H}_{\hbar,V}}^{s})}\lesssim 	\|g_{\varepsilon}\|_{C([0,T];\mathrm{H}_{\mathcal{H}_{\hbar,V}}^{2+s})},
	\end{equation*}
	\sloppy with the constant independent of $\hbar$ and $t\in[0,T]$.
	Since
	$g_{\varepsilon} \rightarrow 0$ in $C([0, T] ; \mathrm{H}_{\mathcal{H}_{\hbar,V}}^{2+s})$, therefore, we have
	\begin{equation*}
		w_{\varepsilon} \to 0 \text{ in }C([0, T]; \mathrm{H}_{\mathcal{H}_{\hbar,V}}^{2+s}),\quad \varepsilon\to 0,
	\end{equation*}
	i.e.,
	\begin{equation*}
		u_{\varepsilon} \to u \text{ in }C([0, T]; \mathrm{H}_{\mathcal{H}_{\hbar,V}}^{2+s}) ,\quad \varepsilon\to 0.
	\end{equation*} 
	This concludes the proof.
\end{proof}
\section{Semi-classical limit $\hbar\to 0$}\label{sec:semiclassic}
In this section,  we mainly compare the classical (resp. very weak) solutions of (\ref{mainpde}) on $\hbar \mathbb{Z}^n$  with the classical (resp. very weak) solutions of (\ref{Eucledian}) on $\mathbb{R}^n$  as $\hbar\to 0$. We start with the proof of  Theorem \ref{semlimit}, the semi-classical limit theorems for the classical solution. In this proof, we have used techniques used in the proof of \cite[Theorem 2.9]{Sch:arxiv}.
\begin{proof}[Proof of Theorem \ref{semlimit}]
	Let $u$ and $v$ solve the Cauchy problems
	\begin{equation}\label{CP1}
		\left\{\begin{array}{l}
			\mathbb{D}_{(g)} u(t, k)+a(t)\mathcal{H}_{\hbar,V} u(t, k)=f(t, k), \quad(t, k) \in(0, T] \times \hbar \mathbb{Z}^n, \\
			u(0, k)=u_0(k), \quad k \in \hbar \mathbb{Z}^n,
		\end{array}\right.
	\end{equation}
	and
	\begin{equation}\label{CP2}
		\left\{\begin{array}{l}
			\mathbb{D}_{(g)} v(t, x)+a(t)\mathcal{H}_{V} v(t, x)=f(t, x), \quad(t, x) \in(0, T] \times   \mathbb{R}^n, \\
			u(0, x)=u_0(x), \quad x \in   \mathbb{R}^{n},
		\end{array}\right.
	\end{equation}
	respectively. Using the above Cauchy problems, denoting $w:=u-v$,  we get
	\begin{equation}\label{CPF}
		\left\{
		\begin{array}{ll}
			\mathbb{D}_{(g)}w(t,k)+a(t)\mathcal{H}_{\hbar,V}w(t,k) =a(t)\left(\mathcal{H}_{V}-\mathcal{H}_{\hbar,V}\right)v(t,k),\quad  t\in(0,T],\\
			w(0,k)=0,\quad k\in\hbar\mathbb{Z}^{n}.\\
		\end{array}		\right.
	\end{equation}
	Applying Theorem \ref{mthm1} for the above Cauchy problem and using estimate \eqref{reqesttt} along with the fact that $w(0,k)\equiv0$, we get   
	\begin{multline}\label{cgt1}
		\|w(t,\cdot)\|_{\mathrm{H}_{\mathcal{H}_{\hbar,V}}^{2+s}}+\|	\mathbb{D}_{(g)} w(t,\cdot)\|_{\mathrm{H}_{\mathcal{H}_{\hbar,V}}^{2+s}}
		\leq C_{a}\|a\left(\mathcal{H}_{V}-\mathcal{H}_{\hbar,V}\right)v\|_{C([0,T];\mathrm{H}_{\mathcal{H}_{\hbar,V}}^{2+s})}\\
		\leq C_{a}\|a\|_{C([0,T])}\|\left(\mathcal{H}_{V}-\mathcal{H}_{\hbar,V}\right)v\|_{C([0,T];\mathrm{H}_{\mathcal{H}_{\hbar,V}}^{2+s})},
	\end{multline}
	where the constant $C_{a}$ is given by \eqref{cnstthm}.
	
	Now we will estimate the term $\left\|\left(\mathcal{H}_{V}-\mathcal{H}_{\hbar,V}\right)v\right\|_{C([0,T];\mathrm{H}_{\mathcal{H}_{\hbar,V}}^{2+s})}$. Let $\phi\in C^{4}(\mathbb{R}^{n})$, then using the Lagrange's form of the remainder and Taylor's theorem, we obtain 
	\begin{equation}\label{tylr}
		\phi(\xi+\mathbf{h})=\sum_{|\alpha|\leq 3} \frac{\partial^{\alpha} \phi(\xi)}{\alpha !} \mathbf{h}^{\alpha}+\sum_{|\alpha|=4} \frac{\partial^{\alpha} \phi(\xi+\theta_{\xi}\mathbf{h})}{\alpha !} \mathbf{h}^{\alpha},
	\end{equation} 
	for some $\theta_{\xi}\in (0,1)$ depending on $\xi$. Utilizing the above relation for   $\mathbf{h}=v_{j}$ and $\mathbf{h}=-v_{j}$ for $j=1,2,\dots,n$, we deduce that
	\begin{equation}\label{d3}
		\sum\limits_{j=1}^{n}\delta_{\xi_{j}}^{2}\phi(\xi)=\sum\limits_{j=1}^{n}\phi^{(2v_{j})}(\xi)\\+\frac{1}{4!}\sum\limits_{j=1}^{n}\left(\phi^{(4v_{j})}(\xi+\theta_{j,\xi} v_{j})+\phi^{(4v_{j})}(\xi-\tilde{\theta}_{j,\xi} v_{j})\right),
	\end{equation}
	where $ \theta_{j,\xi},\tilde{\theta}_{j,\xi}\in(0,1),$ $\delta_{\xi_{j}}^{2}$ is the 2nd order central difference operator and $v_{j}$ is the $j^{th}$ standard unit vector in $\mathbb{Z}^{n}$. Now if we define a translation operator $E_{\theta_{j}v_{j}}\phi:\mathbb{R}^{n}\to \mathbb{R}$  by $E_{\theta_{j}v_{j}}\phi(\xi):=\phi(\xi-\theta_{j,\xi}v_{j}),$ then the equation \eqref{d3} takes the form
	\begin{equation*}
		\sum\limits_{j=1}^{n}\delta_{\xi_{j}}^{2}\phi(\xi)-\sum\limits_{j=1}^{n}\frac{\partial^{2}}{\partial\xi_{j}^{2}}\phi(\xi)=\frac{1}{4!}\sum\limits_{j=1}^{n}\left(E_{-\theta_{j}v_{j}}\phi^{(4v_{j})}(\xi)+E_{\tilde{\theta}_{j}v_{j}}\phi^{(4v_{j})}(\xi)\right).
	\end{equation*}
	Now we extend this to $\hbar\mathbb{Z}^{n}$. Consider a function  $\phi_{\hbar}:\mathbb{R}^{n}\to \mathbb{R}$  defined by
	$\phi_{\hbar}(\xi):=\phi(\hbar\xi)$. Clearly  $\phi_{\hbar}\in C^{4}(\mathbb{R}^{n})$, if we take $\phi\in C^{4}(\mathbb{R}^{n})$. Now we have
	\begin{eqnarray}\label{eqhzn}
		\mathcal{L}_{1}\phi_{\hbar}(\xi)-\mathcal{L}\phi_{\hbar}(\xi)=\frac{1}{4!}\sum\limits_{j=1}^{n}\left(E_{-\theta_{j}v_{j}}\phi_{\hbar}^{(4v_{j})}(\xi)+E_{\tilde{\theta}_{j}v_{j}}\phi_{\hbar}^{(4v_{j})}(\xi)\right),
	\end{eqnarray}
	where $\mathcal{L}$ is the Laplacian on $\mathbb{R}^{n}$ and $\mathcal{L}_{1}$ is the discrete difference Laplacian on $\mathbb{Z}^{n}$.  One can quickly notice that 
	\begin{multline}
		E_{-\theta_{j}v_{j}}\phi_{\hbar}^{(4v_{j})}(\xi)=\phi_{\hbar}^{(4v_{j})}(\xi+\theta_{j,\xi} v_{j})=\hbar^{4}\phi^{(4v_{j})}(\hbar\xi+\hbar\theta_{j,\xi} v_{j})=\\\hbar^{4}E_{-\hbar\theta_{j}v_{j}}\phi^{(4 v_{j})}(\hbar\xi).
	\end{multline}
	Therefore, the equality \eqref{eqhzn} becomes
	\begin{eqnarray}\label{lapdif}
		\left(\mathcal{L}_{\hbar}-\hbar^{2}\mathcal{L}\right)\phi(\hbar\xi)=\frac{\hbar^{4}}{4!}\sum\limits_{j=1}^{n}\left(E_{-\hbar\theta_{j}v_{j}}\phi^{(4v_{j})}(\hbar\xi)+E_{\hbar\tilde{\theta}_{j}v_{j}}\phi^{(4v_{j})}(\hbar\xi)\right).
	\end{eqnarray}
	Combining \eqref{dhamil}, \eqref{eucd} and \eqref{lapdif}, we get 
	\begin{multline}
		\left(\mathcal{H}_{V}-\mathcal{H}_{\hbar,V}\right)\phi(\hbar\xi)=  	\left(\hbar^{-2}\mathcal{L}_{\hbar}-\mathcal{L}\right)\phi(\hbar\xi)=\frac{\hbar^{2}}{4!}\sum\limits_{j=1}^{n}\left(E_{-\hbar\theta_{j}v_{j}}\phi^{(4v_{j})}(\hbar\xi)\right.+\\
		\left.E_{\hbar\tilde{\theta}_{j}v_{j}}\phi^{(4v_{j})}(\hbar\xi)\right).
	\end{multline}
	Hence, it follows that
	\begin{equation}\label{EQ:convh}
		\left\|\left(\mathcal{H}_{V}-\mathcal{H}_{\hbar,V}\right)\phi\right\|^{2}_{\mathrm{H}_{\mathcal{H}_{\hbar,V}}^{2+s}}\lesssim\hbar^{4}\max_{1\leq j\leq n}\left(\left\|E_{-\hbar\theta_{j}v_{j}}\phi^{(4v_{j})}\right\|^{2}_{\mathrm{H}_{\mathcal{H}_{\hbar,V}}^{2+s}}
		+\left\|E_{\hbar\tilde{\theta}_{j}v_{j}}\phi^{(4v_{j})}\right\|^{2}_{\mathrm{H}_{\mathcal{H}_{\hbar,V}}^{2+s}}\right).
	\end{equation}
	Since $V$ is a non-negative polynomial, therefore, the inclusion \eqref{semb} holds, that is,
	\begin{equation*}
		\mathrm{H}_{\mathcal{H}_{V}}^{s}(\mathbb{R}^{n})\subseteq \mathrm{H}^{s}(\mathbb{R}^{n}),\quad  s>0.
	\end{equation*}
	Also, we have the usual Sobolev embeddding:
	\begin{equation*}
		s>k+\frac{n}{2} \Longrightarrow 	\mathrm{H}^{s}(\mathbb{R}^{n}) \subseteq C^k\left(\mathbb{R}^n\right),
	\end{equation*}
	see \cite[Excercise 2.6.17]{Ruz:book}. Combining the above embeddings, we get
	\begin{equation}\label{sobemb}
		s>k+\frac{n}{2} \Longrightarrow 	\mathrm{H}_{\mathcal{H}_{V}}^{s}(\mathbb{R}^{n}) \subseteq C^k\left(\mathbb{R}^n\right).
	\end{equation}
	Since $u_{0} \in \mathrm{H}_{\mathcal{H}_{V}}^{2+s} $ with $s>2+\frac{n}{2}$, therefore,  from Theorem \ref{eucclass}, it follows that the classical solution $v\in C([0,T];\mathrm{H}_{\mathcal{H}_{V}}^{2+s})$ with $s>2+\frac{n}{2}$. Using the embedding \eqref{sobemb}, we get $v\in C^{4}(\mathbb{R}^{n})$ and using the hypothesis $u_{0}^{(4v_{j})} \in \mathrm{H}_{\mathcal{H}_{V}}^{2+s} $ for all $j=1,\dots,n$, we deduce that
	\begin{equation}\label{v4der}
		v^{(4v_{j})}(t,\cdot)\in \mathrm{H}_{\mathcal{H}_{V}}^{2+s}(\mathbb{R}^n), \quad \text{for all }t\in[0,T].
	\end{equation} 
	Now from  \eqref{EQ:convh}, it follows that
	\begin{multline}\label{EQ:cht}
		\left\|\left(\mathcal{H}_{V}-\mathcal{H}_{\hbar,V}\right)v\right\|^{2}_{_{C([0,T];\mathrm{H}_{\mathcal{H}_{\hbar,V}}^{2+s})}}\lesssim\\\hbar^{4}\max_{1\leq j\leq n}\left(\left\|E_{-\hbar\theta_{j}v_{j}}v^{(4v_{j})}\right\|^{2}_{_{C([0,T];\mathrm{H}_{\mathcal{H}_{\hbar,V}}^{2+s})}}+\left\|E_{\hbar\tilde{\theta}_{j}v_{j}}v^{(4v_{j})}\right\|^{2}_{_{C([0,T];\mathrm{H}_{\mathcal{H}_{\hbar,V}}^{2+s})}}\right).
	\end{multline}
	Using \eqref{cgt1}, \eqref{v4der} and \eqref{EQ:cht}, we get $\|w(t,\cdot)\|_{\mathrm{H}_{\mathcal{H}_{\hbar,V}}^{2+s}}+\|\mathbb{D}_{(g)}w(t,\cdot)\|_{\mathrm{H}_{\mathcal{H}_{\hbar,V}}^{2+s}}\to 0$ as $\hbar\to 0.$ Hence $\|w(t,\cdot)\|_{\mathrm{H}_{\mathcal{H}_{\hbar,V}}^{2+s}} \to 0$ and $\|\mathbb{D}_{(g)}w(t,\cdot)\|_{\mathrm{H}_{\mathcal{H}_{\hbar,V}}^{2+s}} \to 0$ as $\hbar \to 0$. This concludes the proof of Theorem \ref{semlimit}.	
\end{proof}
%\begin{proof}[Proof of Theorem \ref{vvyksemlimit}]
The proof of Theorem \ref{vvyksemlimit} will be almost identical to the proof of Theorem \ref{semlimit} with minimal modifications.
%\end{proof}
\section{Remarks}\label{remark}
A few final comments regarding the outcomes of classical and very weak solution in the Euclidean framework are as follows:
\begin{enumerate}
	\item The use of $\mathcal{H}_{V}$-Fourier transform instead of $\mathcal{H}_{\hbar,V}$-Fourier transfrom and following the technique 
	of proof of Theorem \ref{mthm1} will give the proof of Theorem \ref{eucclass}. %A particular example of such transform can be traced in \cite{ldn1} for the Landau Hamiltonian operator.
	\item The notion of unique very weak solution in Theorem \ref{Euwk} can be adapted verbatim by using the Sobolev norm $\mathrm{H}_{\mathcal{H}_{V}}^{s}$ instead of $\mathrm{H}_{\mathcal{H}_{\hbar,V}}^{s}$ in Definition \ref{vwkdef} and Definition \ref{uniquedef}. Consequently, the proof of Theorem \ref{Euwk} will follow along the lines of the proof of Theorem \ref{ext}.
	\item Suppose $\mathcal{M}$ is a positive self-adjoint operator on a Hilbert space $\mathcal{H}$ with a discrete spectrum. Then Theorem \ref{eucclass} can also be extended for the   following Cauchy problem : \begin{equation*}
		\left\{\begin{array}{l}
			\mathbb{D}_{(g)} u(t)+a(t)\mathcal{M} u(t)=f(t), \quad t \in(0, T], \\
			u(0)=u_0   \in \mathcal{H},
		\end{array}\right.
	\end{equation*}
	which is a generalization of \cite{2023niyaz2}.
\end{enumerate}
\bibliographystyle{alphaabbr}
\bibliography{time-fractional}

\end{document}